\documentclass[11pt]{elsarticle}

\usepackage{psfrag,graphicx,booktabs}
\usepackage{amsmath,amssymb,amsthm}
\usepackage{hyperref}

\newtheorem{theorem}{Theorem}[section]

\newtheorem{proposition}[theorem]{Proposition}
\newtheorem{lemma}[theorem]{Lemma}
\newtheorem{example}[theorem]{Example}
\newtheorem{remark}[theorem]{Remark}

\newcounter{appendixsection}
\makeatletter
\gdef\theappendixsection{\@Alph\c@section}
\makeatother
\newtheorem{lemma_appendix}{Lemma}[appendixsection]

\theoremstyle{definition}
\newtheorem{definition}[theorem]{Definition}
\newtheorem{assumption}[theorem]{Assumption}

\newcommand{\X}{X}
\newcommand{\Y}{Y}
\newcommand{\rf}{\mathcal R}
\newcommand{\sif}{\mathcal S}
\newcommand{\sset}{\Sigma}
\DeclareMathOperator{\Fo}{\mathbf F}
\newcommand{\eps}{\varepsilon}
\newcommand{\respar}{\beta}
\newcommand{\req}[1]{\eqref{eq:#1}}

\DeclareMathOperator{\domain}{dom}
\DeclareMathOperator{\range}{ran}

\makeatother
\newcommand{\cl}[1][\tau]{\operatorname{#1-cl}}

\DeclareMathOperator*{\argmin}{arg\,min}

\DeclareMathOperator{\supp}{supp}
\DeclareMathOperator*{\Limsupt}{\tau-Lim\,sup}
\DeclareMathOperator*{\Liminft}{\tau-Lim\,inf}
\DeclareMathOperator*{\Limt}{\tau-Lim}

\DeclareMathOperator*{\Limsuptr}{\tau_{\rf}-Lim\,sup}

\DeclareMathOperator*{\Limtr}{\tau_{\rf}-Lim}

\DeclareMathOperator*{\Limsuptp}{\tau_p-Lim\,sup}

\newcommand{\field}[1]{\ensuremath{\mathbb{#1}}}
\newcommand{\R}{\field{R}}
\newcommand{\N}{\field{N}}

\newcommand{\inner}[2]{\left\langle#1,#2\right\rangle}
\newcommand{\minner}[2]{\bigl\langle#1,#2\bigr\rangle}
\newcommand{\abs}[1]{\left\lvert#1\right\rvert}
\newcommand{\norm}[1]{{\left\lVert#1\right\rVert}}
\newcommand{\mnorm}[1]{\bigl\lVert#1\bigr\rVert}
\newcommand{\snorm}[1]{\lVert#1\rVert}
\newcommand{\set}[1]{\left\{#1\right\}}
\newcommand{\kl}[1]{\left(#1\right)}
\newcommand\enorm{\left\Vert\,\cdot\,\right\Vert}

\begin{document}

\begin{frontmatter}

\title{The Residual Method for Regularizing Ill-Posed Problems}

\author[label1]{Markus Grasmair}
\ead{markus.grasmair@univie.ac.at}
\author[label2]{Markus Haltmeier}
\ead{markus.haltmeier@mpibpc.mpg.de}
\author[label1,label3]{Otmar Scherzer}
\ead{otmar.scherzer@univie.ac.at}

\address[label1]{Computational Science Center, University of Vienna, Nordbergstra{\ss}e 15, Vienna, Austria}
\address[label2]{Statistical Inverse Problems Group, Max Planck Institute for Biophysical Chemistry, G\"ottingen, Germany}
\address[label3]{Johann Radon Institute for Computational and Applied Mathematics, Austrian Academy of Sciences, Austria}

\begin{abstract}
  Although the \emph{residual method}, or \emph{constrained regularization},
  is frequently used in applications, a detailed study of its properties is still missing.
  This sharply contrasts the progress of the theory of Tikhonov regularization,
  where a series of new results for regularization in Banach spaces
  has been published in the recent years.
  The present paper intends to bridge the gap between the existing theories
  as far as possible.
  We develop a stability and convergence theory for the residual method in general topological spaces.
  In addition, we prove convergence rates in terms of (generalized)
  Bregman distances, which can also be applied to non-convex regularization
  functionals.

  We provide three examples that show the applicability of our theory.
  The first example is the regularized solution of linear operator equations
  on $L^p$-spaces, where we show that the results of Tikhonov regularization
  generalize unchanged to the residual method.
  As a second example, we consider the problem of density estimation
  from a finite number of sampling points,
  using the Wasserstein distance as a fidelity term and an entropy measure
  as regularization term.
  It is shown that the densities obtained in this way depend continuously
  on the location of the sampled points and that the underlying density can
  be recovered as the number of sampling points tends to infinity.
  Finally, we apply our theory to compressed sensing.
  Here, we show the well-posedness of the method and derive convergence rates
  both for convex and non-convex regularization under rather weak conditions.
\end{abstract}

\begin{keyword}
Ill-posed problems \sep 
Regularization \sep
Residual Method \sep 
Sparsity \sep 
Stability \sep 
Convergence Rates.

\MSC
65J20\sep 
47J06\sep 
49J27.
\end{keyword}

\end{frontmatter}

\section{Introduction}

We study the solution of ill-posed operator equations
\begin{equation}\label{eq:F}
  \Fo(x) = y\,,
\end{equation}
where $\Fo\colon \X \to \Y$ is an operator between the topological spaces $\X$ and $\Y$,
and $y \in \Y$ are given, noisy data,
which are assumed to be close to some unknown, noise-free data
$y^\dagger \in \range(\Fo)$.
If the operator $\Fo$ is not continuously invertible,
then~\req{F} may not have a solution and,
if a solution exists, arbitrarily small perturbations
of the data may lead to unacceptable results.

If $\Y$ is a Banach space and the given data are known to satisfy an estimate
$\norm{y^\dagger-y} \le \respar$, one strategy for defining an approximate
solution of~\req{F} is to solve the \emph{constrained minimization problem}
\begin{equation}\label{eq:constrained}
  \rf(x) \to \min
  \qquad \text{ subject to }
   \quad \norm{\Fo(x)-y} \le \respar\;.
\end{equation}
Here, the \emph{regularization term} $\rf\colon \X \to [0,+\infty]$
is intended to enforce certain regularity properties of the approximate solution
and to stabilize the process of solving~\req{F}.
In~\cite{IvaVasTan02,Tan97}, this strategy is called the \emph{residual method}.
It is closely related to \emph{Tikhonov regularization}
\begin{equation}\label{eq:tik}
  {\cal T}(x) := \norm{\Fo(x)-y}^2 + \alpha\rf(x) \to \min \,,
\end{equation}
where $\alpha > 0$ is a regularization parameter. In the case that the operator $\Fo$ is linear
and $\rf$ is convex, \req{constrained} and \req{tik} are basically
equivalent, if $\alpha$ is chosen according to Morozov's discrepancy principle (see  \cite[Chap. 3]{IvaVasTan02}).

While the theory of Tikhonov regularization  has received much attention in the literature
(see for instance \cite{AcaVog94,DauDefDem04,EngHanNeu96,EngKunNeu89,Gro84,HofYam05,Mor93,SchEngKun93,
SeiVog89,TikArs77,Vas06}),
the same cannot be said about the residual method.
The existing results are mainly concerned with the existence theory of~\req{constrained} and with the question
of convergence, which asks whether solutions of~\req{constrained}
converge to a solution of~\req{F} as $\norm{y - y^\dagger} \le \respar \to 0$.
These problems have been treated in very general settings in~\cite{Iva69,Sei81}
(see also~\cite{Hei08,Tan95,Tan97}).
Convergence rates have been derived in~\cite{BurOsh04} for linear equations in Hilbert spaces and later generalized
in~\cite{Hei08} to non-linear equations in Banach spaces.
Convergence rates have also been derived 
in~\cite{Can08,CanRomTao06b,GraHalSch11} for the reconstruction of sparse sequences.

The problem of stability, however, that is,
continuous dependence of the solution of~\req{constrained} on the input
data $y$ and the presumed noise level $\respar$,
has been hardly considered at all.
One reason for the lack of results is that, in contrast to Tikhonov regularization, stability
simply does not hold for general non-linear operator equations.
But even for the linear case, where we indeed prove stability,
so far stability theorems are non-existent in the literature.
Though some results have been derived in~\cite{Hei08},
they only cover a very weak form of stability,
which states that the solutions of~\req{constrained} with perturbed data
stay close to the solution with unperturbed data,
if one additionally increases the regularization parameter $\respar$ in the
perturbed problem by a sufficient amount.

\medskip
The present paper tries to generalize
the existent theory on the residual method as far as possible.
We assume that $\X$ and $\Y$ are mere topological spaces and consider the minimization
of $\rf(x)$ subject to the constraint $\sif\left(\Fo(x),y\right) \le \respar$.
Here $\sif$ is some distance like functional taking over the role of the norm in~\req{constrained}.
In addition, we discuss the case where the operator $\Fo$ is not known exactly.
This subsumes errors due to the modeling process as well as discretizations of the problem necessary for its numerical solution.
We provide different criteria that ensure
stability (Lemma \ref{le:stability}, Theorem \ref{thm:stability:2} and Proposition~\ref{pr:linear_stab})
and convergence (Theorem~\ref{thm:convergence} and Proposition~\ref{pr:linear_stab})
of the residual method.
In particular, our conditions also include certain non-linear operators
(see Example~\ref{ex:stab_nonlinear}).

Section~\ref{se:rates} is concerned with the derivation of convergence rates,
i.e.,  quantitative estimates between solutions of \req{constrained} and
the exact data $y^\dag$.
Using notions of abstract convexity, we define a generalized Bregman distance that allows us
to state and prove rates on arbitrary topological spaces (see Theorem~\ref{thm:rates}).
In Section~\ref{sec:sparsity} we apply our general results to the case of sparse $\ell^p$-regularization
with $p \in (0,2)$.
We prove the well-posedness of the method and derive convergence rates
with respect to the norm in a fairly general setting.
In the case of convex regularization, that is, $p \ge 1$,
we derive a convergence rate of order $\mathcal{O}\bigl(\respar^{1/p}\bigr)$.
In the non-convex case $0 < p < 1$, we show that the rate $\mathcal{O}(\respar)$ holds.

\section{Definitions and Mathematical Preliminaries}\label{sec:defs}

Throughout the paper, $\X$ and $\Y$ denote sets.
Moreover, $\rf \colon \X \to [0,+\infty]$ is a functional on $\X$,
and  $\sif\colon \Y \times \Y \to [0,+\infty]$ is a
functional on $\Y\times \Y$ such that $\sif(y,z) = 0$ if and only if $y = z$.

\subsection{The Residual Method}
For given mapping $\Fo\colon \X \to \Y$, given data $y \in \Y$, and fixed  parameter
$\respar \ge 0$, we consider the constrained minimization problem
\begin{equation}\label{eq:min_prob}
  \rf(x) \to \min
  \qquad\text{ subject to } \quad
  \sif(\Fo(x),y) \le \respar\;.
\end{equation}
For the analysis of the residual method \req{min_prob} it is convenient to
introduce the following notation.

The \emph{feasible set} $\Phi(\Fo,y,\respar)$, the \emph{value} $v(\Fo,y,\respar)$,
and the \emph{set of solutions} $\sset(\Fo,y,\respar) $ of \req{min_prob} are defined by
\begin{align*}
    \Phi(\Fo,y,\respar)   &:= \set{x \in \X: \sif(\Fo(x),y)\le \respar}\,,\\
    v(\Fo,y,\respar)      &:= \inf\set{\rf(x):x \in \Phi(\Fo,y,\respar)} \,, \\
   \sset(\Fo,y,\respar)  &:= \set{x \in \Phi(\Fo,y,\respar): \rf(x) = v(\Fo,y,\respar)}  \,.
\end{align*}
In particular, $\Phi(\Fo,y,0)$ consist of all solutions of the equation $\Fo(x) = y$.
The elements of $\sset(\Fo,y,0)$ are therefore referred to as
$\rf$-\emph{minimizing solutions}  of $\Fo(x) = y$.

In addition,  for $t \geq 0$, we set
\begin{equation}\label{eq:fr}
	\Phi_\rf(\Fo,y,\respar,t) := \Phi(\Fo,y,\respar) \cap \set{x \in \X: \rf(x) \le t} \,.
\end{equation}
An immediate consequence of the above definitions is the identity
\begin{equation}\label{eq:Sigma_min}
    \sset(\Fo,y,\respar) = \Phi_\rf \left(\Fo,y,\respar, v(\Fo,y,\respar)\right) \,.
\end{equation}

\begin{remark}
We do not assume a--priori that a solution of the minimization problem \req{min_prob} exists.
Only in the next section shall we deduce the existence of solutions under a
compactness assumption on the sets $\Phi_\rf(\Fo,y,\respar,t)$, see Theorem \ref{th:existence}.
\end{remark}

\begin{lemma}\label{le:Phi_prop}
The sets  $\Phi_{\rf}(\Fo,y,\respar,t)$ defined in \req{fr} satisfy
  \begin{equation}\label{eq:Sigma_order}
  \Phi_\rf(\Fo,y,\respar,t) \subset \Phi_\rf(\Fo,y,\respar+\gamma,t+\eps)
  \end{equation}
for every $\gamma, \eps \ge 0$, and
 \begin{equation}\label{eq:Sigma_intersect}
  \Phi_{\rf}(\Fo,y,\respar,t) =
  \bigcap_{\gamma,\eps > 0} \Phi_{\cal R}(\Fo,y,\respar+\gamma,t+\eps)\;.
  \end{equation}
\end{lemma}

\begin{proof}
  The inclusion \req{Sigma_order} follows immediately from the definition of $\Phi_{\cal R}$.
  For the proof of~\req{Sigma_intersect} note that   $x \in \bigcap_{\gamma,\eps > 0} \Phi_{\cal R}(\Fo,y,\respar+\gamma,t+\eps)$
  if and only if $\sif (\Fo(x),y) \le \respar+\gamma$ for all $\gamma > 0$
  and $\rf(x) \le t+\eps$ for all $\eps > 0$.
  This, however, is the case if and  only if  $\sif(\Fo(x),y) \le \respar$
  and $\rf(x) \le t$, which means that $x \in \Phi_{\rf}(\Fo,y,\respar,t)$.
\end{proof}

Further properties of the value $v$ and the sets $\Phi$, $\Phi_{\rf}$ and $\sset$
are summarized in \ref{sec:appendix}.

\subsection{Convergence of Sets of Solutions}

In the next section we study convergence and stability of the residual method,
that is, the behavior of the set of solutions $\sset(\Fo_k,y_k,\respar_k)$ for $\respar_k \to \respar$,
$y_k \to y$, and $\Fo_k \to \Fo$.
In~\cite{EngHanNeu96,SchGraGroHalLen09}, where convergence and stability of
Tikhonov regularization have been investigated, the stability results are of the following form:
For every sequence $(y_k)_{k\in\N} \to y$ and every sequence of minimizers
$x_k \in \argmin \bigl\{ \norm{\Fo(x)-y_k}^2 + \alpha\rf(x) \bigr\}$ there exists a subsequence of
$(x_k)_{k\in\N}$ that converges to a minimizer of $\norm{\Fo(x)-y}^2 + \alpha\rf(x)$.
In this paper we prove similar results for the residual method but with a different notation
using a type of convergence of sets (see, for example, \cite[\S 29]{Kur66}).

\begin{definition}\label{de:Limsup}
  Let $\tau$ be a topology on $\X$ and let $(\sset_k)_{k\in\N}$  be a sequence
  of subsets of $\X$.
  \begin{enumerate}
  \item[(a)]
  The \emph{upper limit} of $(\sset_k)_{k\in\N}$ is defined as
    \[
    \Limsupt_{k\to \infty} \sset_k := \bigcap_{k\in\N} \biggl(\cl \bigcup_{k'\ge k} \sset_{k'}\biggr)\;,
    \]
    where $\cl$ denotes the closure with respect to $\tau$.

  \item[(b)]
    An element $x \in \X$ is contained in the lower limit
    of the sequence $(\sset_k)_{k\in\N}$, in short
    \[
    x \in \Liminft_{k\to \infty} \sset_k \,,
    \]
    if for every neighborhood $N$ of $x$ there exists $k_0\in\N$ such that
    $N\cap \sset_{k} \neq \emptyset$ for every $k\ge k_0$.

  \item[(c)]
  If the lower limit and the upper limit of $(\sset_k)_{k\in\N}$ coincide,
  we define
    \[\Limt_{k\to\infty} \sset_k := \Liminft_{k\to \infty} \sset_k = \Limsupt_{k\to \infty} \sset_k \]
    as the \emph{limit} of the sequence $(\sset_k)_{k\in\N}$.
\end{enumerate}
\end{definition}

\begin{remark}\label{re:limsup}
  As a direct consequence of Definition~\ref{de:Limsup}, an element $x$ is contained in the
  upper limit $\Limsupt_{k\to \infty} \sset_k$, if and only
  if for every neighborhood $N$ of $x$ and every $k_0 \in \N$ there exists $k \ge k_0$ 
  with $N\cap \sset_{k} \neq \emptyset$.

  If $\X$ satisfies the first axiom of countability,
  then $x \in \Limsupt_{k\to \infty} \sset_k$, if and only if there exists
  a subsequence $(\sset_{k_j})_{j\in\N}$ of $(\sset_{k})_{k\in\N}$ and a sequence
  of elements $x_j \in \sset_{k_j}$ such that $x_j \to_\tau x$ (see~\cite[\S 29.IV]{Kur66}).
  Note that in particular every metric space satisfies the first axiom of countability.
\end{remark}

The following proposition clarifies the relation between the stability and convergence results
in~\cite{EngHanNeu96,SchGraGroHalLen09} and the results in the present paper.

\begin{proposition}
  Let $(\sset_k)_{k \in \N}$ be a sequence of nonempty subsets of $\X$, and assume that there exists
  a compact set $K$ such that $\sset_k \subset K$ for all $k \in \N$.
  Then $\Limsupt_{k\to \infty} \sset_k$ is non-empty.

  If, in addition, $\X$ satisfies the first axiom of countability,
  then every sequence of elements $x_k \in \sset_k$ has a subsequence converging
  to some element $x \in \Limsupt_{k\to \infty} \sset_k$.
\end{proposition}

\begin{proof}
  By assumption, the sets $S_k := \cl \bigcup_{k'\ge k}\sset_k$ form a decreasing
  family of non-empty, compact sets.
  Thus also their intersection $\bigcap_{k\in\N} S_k = \Limsupt_{k\to \infty} \sset_k$
  is non-empty (see~\cite[Thm.~5.1]{Kel55}).

  Now assume that $\X$ satisfies the first axiom of countability.
  Then in particular every compact set is sequentially compact (see~\cite[Thm.~5.5]{Kel55}).
  Let now $x_k \in \sset_k$ for every $k\in\N$.
  Then $(x_k)_{k\in\N}$ is a sequence in the compact set $K$ and therefore
  has a subsequence $(x_{k_j})_{j\in\N}$ converging to some element $x \in K$.
  From Remark~\ref{re:limsup} it follows that $x \in \Limsupt_{k\to\infty} \sset_k$,
  which shows the assertion.
\end{proof}

\subsection{Convergence of the Data}

In addition to the convergence of subsets $\sset_k$ of $\X$,
it is necessary to define a notion of convergence on the set
$\Y$ that is compatible with the distance measure $\sif$.

\begin{definition}\label{de:s_uniform}
The sequence $(y_k)_{k\in\N} \subset \Y$ converges \emph{$\sif$-uniformly} to $y \in \Y$, if
  \[
  \sup\set{\abs{\sif(z,y_k)-\sif(z,y)}:z \in \Y} \to 0\;.
  \]
The sequence of mappings $\Fo_k \colon \X \to \Y$  converges \emph{locally $\sif$-uniformly}
to $\Fo\colon \X \to \Y$, if
  \[
  \sup\set{\abs{\sif(\Fo_k(x),y)-\sif(\Fo(x),y)}: y \in \Y,\, x \in \X,\, \rf(x) \le t} \to 0
  \]
for every $t \ge 0$.
\end{definition}

\begin{remark}
The $\sif$-uniform convergence on $\Y$ is induced by the extended metric
\[ \sif_1(y_1,y_2) := \sup\set{\abs{\sif(z,y_1)-\sif(z,y_2)} : z \in \Y } \,. \]
If the distance measure $\sif$ itself equals a metric, then $\sif_1$ coincides
with $\sif$. Similarly, local $\sif$-uniform convergence of  a sequence of mappings $\Fo_k$
equals the uniform convergence of $\Fo_k$ on $\rf$-bounded sets with respect to the extended metric \[\sif_2(y_1,y_2) := \sup\set{ \abs{\sif(y_1,z)-\sif(y_2,z)}:z \in \Y } \,. \]
\end{remark}

\section{Well-posedness of the Residual Method} \label{sec:wp}

In the following we investigate the existence of minimizers,
and the stability and the convergence of the residual method.
Throughout the whole section we assume that $\tau$ is a topology on $\X$,
$\Fo\colon \X \to \Y$ is a mapping, $y \in \Y$ are given data
and $\respar \ge 0$ is a fixed parameter.

\subsection{Existence}

We first investigate under which conditions $\sset(\Fo,y,\respar)$,
the set of solutions of  \req{min_prob}, is not empty.

\begin{theorem}[Existence]\label{th:existence}
  Assume that $\Phi_\rf(\Fo,y,\respar,t)$ is $\tau$-compact for every $t\geq 0$
  and non-empty for some $t_0 \geq 0$.
  Then Problem~\req{min_prob} has a solution.
\end{theorem}

\begin{proof}
Equation \req{Sigma_min} and Lemma \ref{le:Phi_prop} imply the identity
\[
    \sset(\Fo,y,\respar) = \Phi_\rf\kl{\Fo,y,\respar,v(\Fo,y,\respar)}
    =
    \bigcap_{\eps > 0} \Phi_\rf \kl{\Fo,y,\respar,v(\Fo,y,\respar)+\eps} \,.
\]
Because $\Phi_\rf(\Fo,y,\respar,t_0) \neq \emptyset $, the value of \req{min_prob}
satisfies  $v(\Fo,y,\respar) \leq t_0 < \infty$ and therefore
$\emptyset \neq \Phi_\rf\left(\Fo,y,\respar,v(\Fo,y,\respar)+\eps\right)$ for every $\eps>0$.
Consequently, $\sset(\Fo,y,\respar)$ is the intersection of a decreasing family of non-empty $\tau$-compact
sets and thus non-empty (see~\cite[Thm.~5.1]{Kel55}).
\end{proof}

  Recall that a mapping $\mathcal{F}\colon \X \to [0,+\infty]$ is
  \emph{lower semi-continuous},
  if its lower level sets $\set{x \in \X:\mathcal{F}(x) \le t}$ are closed for every $t \ge 0$.
  Moreover, the mapping $\mathcal{F}$ is \emph{coercive}, if its lower level sets are
  pre-compact, see \cite{Bra02}.
  (In a Banach space  one often  calls a functional coercive,
  if it is unbounded on unbounded sets. The notion used here is equivalent if
  the Banach space is reflexive and $\tau$ is the weak topology.)
  In particular, the mapping  $\mathcal{F}$ is  \emph{lower semi-continuous and coercive}, if  and only
  if its lower level sets are compact.

  \begin{proposition}\label{prop:compact}
  Assume that  $\rf$ and $x \mapsto \sif(\Fo(x),y)$ are lower semi-continuous and
  one of them, or their sum, is coercive.
  Then  $\Phi_\rf(\Fo,y,\respar,t)$ is $\tau$-compact for every $t\geq 0$.
  If additionally $\Phi_\rf(\Fo,y,\respar,t_0)$ is non-empty for some $t_0 \geq 0$, then
  Problem~\req{min_prob} has a solution.
  \end{proposition}

  \begin{proof}
  If $\rf$ and $x \mapsto \sif(\Fo(x), y)$
  are lower semi-continuous and one of them is coercive, then
  \[\Phi_\rf(\Fo,y,\respar,t) = \set{x: \sif(\Fo(x), y) \leq \respar} \cap \set{x : \rf(x) \le t}\]
  is the intersection of a closed and a $\tau$-compact set and therefore itself $\tau$-compact.
  In case that only the sum $x \mapsto \sif(\Fo(x), y) + \rf(x)$ is coercive, the set
  \begin{multline*}
    \Phi_\rf(\Fo,y,\respar,t) = \set{x: \sif(\Fo(x), y) \leq \respar} \cap \set{x : \rf(x) \le t}
    \\ \subset
    \{x: \sif(\Fo(x), y) + \rf(x) \leq \respar +t\}
     \end{multline*}
  is a closed set contained in a $\tau$-compact set and therefore
  again $\tau$-compact.
  \end{proof}

  The lower semi-continuity of  $x \mapsto \sif\left(\Fo(x),y\right)$ certainly holds
  if $\Fo$ is continuous and $\sif$ is lower semi-continuous with respect to
  the first component (for some given topology on $\Y$).
  It is, however, also possible to obtain lower semi-continuity,
  if $\Fo$ is not continuous but the functional $\sif$ satisfies
  a stronger condition:

  \begin{proposition}\label{le:Slsc}
  Let $\tau'$ be a topology on $\Y$ such that
  $z \mapsto \sif(z,y)$ is lower semi-continuous
  and coercive, and assume that $\Fo \colon \X \to \Y$ has a closed graph.
  Then the functional $x \mapsto \sif(\Fo(x),y)$ is lower semi-continuous.
  \end{proposition}

\begin{proof}
Because $\Fo$ has a closed graph, the pre-image under $\Fo$
of every compact set is closed (see~\cite[Thm.~4]{Iva69}).
This shows that
\[
  \set{x \in \X:\sif(\Fo(x),y) \le \respar} = \Fo^{-1}\left(\set{z \in \Y:\sif(z,y) \le \respar}\right)
\]
is closed  for every $\respar$, that is, the mapping
$x \mapsto \sif(\Fo(x),y)$ is lower semi-continuous.
\end{proof}

\subsection{Stability}

Stability is concerned with the continuous dependence of the solutions
of~\req{min_prob} of the input data, that is, the element $y$,
the parameter $\respar$, and, possibly, the operator $\Fo$.
Given sequences $\respar_k \to \respar$, $y_k \to y$, and $\Fo_k \to \Fo$,
we ask whether the sequence of sets $\sset(\Fo_k,y_k,\respar_k)$ converges
to $\sset(\Fo,y,\respar)$.
As already indicated in Section~\ref{sec:defs}, we will make use of
the upper convergence of sets introduced in Definition~\ref{de:Limsup}.
The topology, however, with respect to which the results are formulated,
is stronger than $\tau$.

\begin{definition}\label{de:tauR}
  The topology $\tau_\rf$ on $\X$ is generated by all
  sets of the form
  $U \cap \{x \in \X:  \rf(x) > s \}$ with   
  $s\in \R$ and $U \in \tau$  and all sets of the form $U \cap \{x \in \X :   \rf(x) <  t\}$ with   
  $t \in \R \cup \set{\infty}$ and $U \in \tau$.
  (Hence $\tau_\rf$ consists of  all unions of finite
  intersections of sets of the form 
  $U \cap \{x \in \X:  \rf(x) > s \}$ or $U \cap \{x \in \X :   \rf(x) <  t\}$.)

Note that a sequence $(x_k)_{k\in\N} \subset \X$ converges to $x$ with respect
to $\tau_\rf$, if and only if $(x_k)_{k\in\N} $ converges to $x$ with respect
to $\tau$ and satisfies $\rf(x_k) \to \rf(x)$ for $k\to \infty$.
\end{definition}

\medskip
For the stability results we make the following assumption:
\begin{assumption}\label{as:stability} \mbox{}
\begin{enumerate}
\item
Let $\respar \geq 0$, let $y \in \Y$,  and let $\Fo \colon \X \to \Y$ be a mapping.

\item
Let $(\respar_k)_{k\in\N}$ be a sequence of nonnegative numbers, let
$(y_k)_{k\in\N}$ be a sequence in $\Y$, and let $(\Fo_k)_{k\in\N}$
be a sequence of mappings $\Fo_k \colon \X \to \Y$.

\item
The sequence $(\respar_k)_{k\in\N}$ converges to $\respar$, the sequence $(y_k)_{k\in\N}$ converges
$\sif$-uniformly to $y$,  and  $(\Fo_k)_{k\in\N}$ converges locally $\sif$-uniformly to $\Fo$.

\item
The sets $\Phi_\rf(\Fo_k,w,\gamma,t)$ and
$\Phi_\rf(\Fo,w,\gamma,t)$ are compact for all
 $w$, $\gamma$, $t$, and $k$.  Moreover, for every   $w$, $\gamma$,  $k$  
 there exist some $t_0$ such that $\Phi_\rf(\Fo_k,w,\gamma,t_0)$ and
$\Phi_\rf(\Fo,w,\gamma,t_0)$ are nonempty.
\end{enumerate}
\end{assumption}

\medskip
The following lemma is the key result to prove stability of the
residual method.

\begin{lemma}\label{le:stability}
  Let Assumption~\ref{as:stability} hold and assume that
  \begin{equation}\label{eq:stability:V}
    \limsup_{k\to\infty} v(\Fo_k,y_k,\respar_k) \le v(\Fo,y,\respar) < \infty\,.
  \end{equation}
  Then,
  \begin{equation}\label{eq:stability:lim}
  \emptyset \neq \Limsuptr_{k\to\infty} \sset(\Fo_k,y_k,\respar_k) \subset \sset(\Fo,y,\respar)\;.
  \end{equation}
  If, additionally, the set $\sset(\Fo,y,\respar)$ consists of a single element $x_\respar$,
  then
  \begin{equation}\label{eq:stability:lim_unique}
    \{x_\respar\} = \Limtr_{k\to\infty} \sset(\Fo_k,y_k,\respar_k)\;.
  \end{equation}
\end{lemma}

\begin{proof}
  In order to simplify the notation, we define
  \[
  \begin{aligned}
    \Phi_k(t) &:= \Phi_\rf(\Fo_k,y_k,\respar_k,t)\,,
    & \Phi(t) &:= \Phi_\rf(\Fo,y,\respar,t)\,,\\
    v_k &:= v(\Fo_k,y_k,\respar_k)\,,
    & v & := v(\Fo,y,\respar)\,,\\
    \sset_k &:= \sset(\Fo_k,y_k,\respar_k)\,,
    & \sset &:= \sset(\Fo,y,\respar)\;.
  \end{aligned}
  \]
  Moreover we define the set $T := \Limsupt_{k\to \infty} \sset_k$.
  Because the topology $\tau_\rf$ is finer than $\tau$,
  it follows that $\Limsuptr_{k\to\infty}  \sset_k \subset T$.
  We proceed by showing that $\emptyset \neq T \subset \sset$
  and $T \subset \Limsuptr_{k\to\infty} \sset_k$,
  which then gives the assertion~\req{stability:lim}.

\smallskip
  The inequality~\req{stability:V} implies that for every $\eps > 0$
  there exists some $k_0 \in \N$ such that $v_k \le v + \eps$ for all $k \ge k_0$.
  Since $\respar_k \to \respar$, we may additionally assume that $\respar_k \le \respar +\eps$.
  Lemma~\ref{le:Phiincl2} implies, after possibly enlarging $k_0$,
  \begin{multline}\label{eq:stability:h1}
    \Phi_k(v_k)
    \subset \Phi_\rf(\Fo_k,y_k,\respar+\eps,v_k)
    \\
    \subset \Phi_\rf(\Fo,y,\respar+2\eps,v_k)
    \subset \Phi_\rf(\Fo,y,\respar+2\eps,v+\eps)
  \end{multline}
  for all $k \ge k_0$.
  Thus,
  \begin{multline}\label{eq:stability:h1a}
    T = \Limsupt_{k\to \infty} \sset_k
    = \bigcap_{k \in \N} \biggl(\cl[\tau] \bigcup_{k'\ge k} \sset_{k'}\biggr)\\
    = \bigcap_{k \ge k_0} \biggl(\cl[\tau] \bigcup_{k'\ge k} \Phi_{k'}(v_{k'})\biggr)
    \subset \Phi_\rf(\Fo,y,\respar+2\eps,v+\eps)\;.
  \end{multline}
  The sets $\cl[\tau]\bigcup_{k'\ge k} \sset_{k'}$ are closed and non-empty and,
  by assumption, the set $\Phi_\rf(\Fo,y,\respar+2\eps,v+\eps)$ is compact.
  Thus $T$ is the intersection of a decreasing family of non-empty compact
  sets and therefore non-empty.
  Moreover, because~\req{stability:h1a} holds for every $\eps > 0$, we have
  \begin{equation}\label{eq:stability:h1b}
    \emptyset\neq T
    \subset \bigcap_{\eps > 0} \Phi_\rf(\Fo,y,\respar+2\eps,v+\eps)
    = \Phi(v)
    = \sset\;.
  \end{equation}

  Next we show the inclusion $T \subset \Limsuptr_{k\to\infty} \sset_k$.
  To that end, we first prove that
  \begin{equation}\label{eq:stability:limV}
    v = \lim_k v_k\;.
  \end{equation}
  Recall that Theorem~\ref{th:existence} implies that $\Phi_k(v_k) = \sset_k$
  is non-empty. Therefore, \req{stability:h1} implies that also
  $\Phi_\rf(\Fo,y,\respar+2\eps,v_k)$ is non-empty, which in turn
  shows that $v_k \ge v(\Fo,y,\respar+2\eps)$ for all $k$ large enough.
  Consequently,
  \begin{equation}\label{eq:stability:h2}
    \liminf_{k\to \infty} v_k \ge v(\Fo,y,\respar+2\eps)
  \end{equation}
  for all $\eps > 0$.
  From Lemma~\ref{le:Hlimit} we obtain that $v = \sup_{\eps > 0} v(\Fo,y,\respar+2\eps)$.
  Together with~\req{stability:h2} and~\req{stability:V} this shows~\req{stability:limV}.

  Let now $x \in T$, let $N$ be a neighborhood of $x$ with respect to $\tau$,
  let $\delta > 0$ and $k_0 \in \N$. Since $T \subset \sset$ 
  (see~\req{stability:h1b}),
  it follows that $\rf(x) = v$. Thus it follows from~\req{stability:limV} that
  there exists $k_1 \ge k_0$ such that
  \[
  \abs{v_k - \rf(x)} < \delta
  \]
  for all $k \ge k_1$.
  In particular,
  \begin{equation}\label{eq:stability:N}
    \sset_k \subset \set{\tilde{x}\in \X:\rf(x)-\delta < \rf(\tilde{x}) < \rf(x)+\delta}
  \end{equation}
  for all $k \ge k_1$.
  Remark~\ref{re:limsup}  implies that there exists $k_2 \ge k_1$ such that
  \begin{equation}\label{eq:stability:V_old}
    N \cap \sset_{k_2} \neq \emptyset\;.
  \end{equation}
  Now recall that the sets
  $N \cap  \set{\tilde{x}\in \X:\rf(x)-\delta < \rf(\tilde{x}) < \rf(x)+\delta}$
  form a basis of neighborhoods of $x$ for the topology $\tau_\rf$.
  Therefore~\req{stability:N}, \req{stability:V_old}, and the characterization
  of the upper limit of sets given in Remark~\ref{re:limsup} imply that
  $x \in \Limsuptr_{k\to\infty} \sset_k$.
  Thus the inclusion~\req{stability:lim} follows.

  If the set $\sset(\Fo,y,\respar)$ consists of a single element $x_\respar$,
  then the first part of the assertion implies that for every subsequence $(k_j)_{j\in\N}$
  we have
  \[
  \Limsuptr_{j\to\infty} \sset(\Fo_{k_j},y_{k_j},\respar_{k_j}) = \{x_\respar\}\;.
  \]
  Thus the assertion follows from Lemma~\ref{le:subseq}.
\end{proof}\smallskip

The crucial condition in Lemma~\ref{le:stability} is the inequality~\req{stability:V}.
Indeed, one can easily construct examples, where this condition fails and
the solution of Problem~\req{min_prob} is unstable, see Example~\ref{ex:stability} below.
What happens in this example is that the upper limit
$\Limsuptr_{k\to\infty} \sset(\Fo,y_k,\respar)$ consists of local minima
of $\rf$ on $\Phi(\Fo,y,\respar)$ that fail to be global minima of $\rf$
restricted to $\Phi(\Fo,y,\respar)$.

\begin{figure}[htb!]
  \psfrag{y}{$y$}
  \psfrag{ypd}{$y+\respar$}
  \psfrag{ymd}{$y-\respar$}
  \psfrag{1}{1}
  \psfrag{0}{0}
  \centering
  \includegraphics[width=0.6\textwidth]{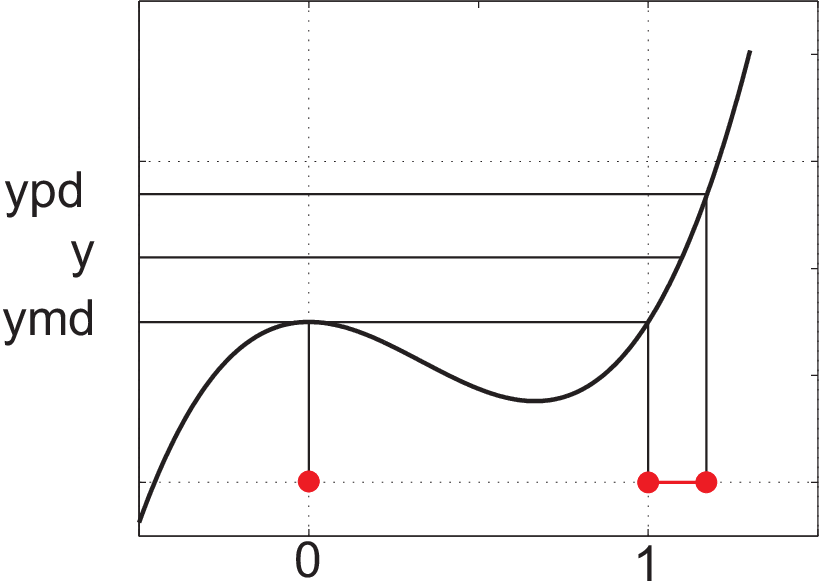}
  \caption{The nonlinear function $\Fo$ from Example \ref{ex:stability}.
    The feasible set $\Phi(\Fo,y,\respar) = \set{x \in \R:  \abs{\Fo(x)-y} \le \respar}$ consists of an interval and the isolated point $\{0\}$.}
\end{figure}

\begin{example}\label{ex:stability}
  Consider the function $\Fo\colon \R \to \R$, $\Fo(x) = x^3 - x^2$,
  and the regularization functional $\rf(x) = x^2$.
  Let $y > 0$ and choose $\respar = y$.
  Then
  \begin{multline} \label{eq:stab_ex}
  \argmin\set{\rf(x):\abs{\Fo(x)-y} \le \respar}
  \\ 
  = \argmin\set{x^2:\abs{x^3-x^2-y} \le y}
  = 0 \;.
  \end{multline}
  Now let $y_k > y$.
  Then
  \[
  \argmin\set{\rf(x):\abs{\Fo(x)-y_k} \le \respar}
  = \argmin\set{x^2:\abs{x^3-x^2-y_k} \le y}
  = x_k\,,
  \]
  where $x_k$ is the unique solution of the equation $\Fo(x) = y_k - y$.
  Thus, if the sequence $(y_k)_{k\in\N}$ converges to $y$ from above,
  we have $x_k > 1$ for all $k$ and $\lim_{k\to \infty} x_k = 1$.
  According to~\req{stab_ex}, however, the solution of the limit problem equals zero.
\end{example}

\medskip
The following two theorems are central results of this paper.
They answer the question to which extent we obtain stability results for the residual
method  similar  to the ones known for  Tikhonov
regularization.

\begin{theorem}[Approximate Stability]\label{thm:ap_stability}
  Let Assumption~\ref{as:stability} hold.
  Then there exists a sequence $\eps_k \to 0$ such that
  \[
  \emptyset \neq \Limsuptr_{k\to\infty} \sset(\Fo_k,y_k,\respar_k+\eps_k) \subset \sset(\Fo,y,\respar)\;.
  \]
\end{theorem}

\begin{proof}
  Define
  \[
  \eps_k := \inf\set{\eps > 0:    \Phi_\rf\left(\Fo,y,\respar,v(\Fo,y,\respar)\right)
    \subset \Phi_\rf\left(\Fo_k,y_k,\respar_k+\eps,v(\Fo,y,\respar)\right)}\;.
  \]
  Lemma~\ref{le:Phiincl2} and the assumption that $\respar_k \to \respar$
  imply that $\eps_k \to 0$.
  Since by assumption
  \[
  \emptyset \neq \sset(\Fo,y,\respar)
  = \Phi_\rf\left(\Fo,y,\respar,v(\Fo,y,\respar)\right)
  \subset \Phi_\rf\left(\Fo_k,y_k,\respar_k+\eps_k,v(\Fo,y,\respar)\right)\,,
  \]
  we obtain that $v(\Fo_k,y_k,\respar_k+\eps_k) \le v(\Fo,y,\respar)$.
  Thus the assertion follows from Lemma~\ref{le:stability}.
\end{proof}

Theorem~\ref{thm:ap_stability}, is a stability result in the same spirit as the
one derived in~\cite{Hei08}.
While it does not assert that, in the general setting described by Assumption~\ref{as:stability},
the residual method is stable in the sense that the solutions depend continuously
on the input data, it does state that the solutions of the perturbed problems
stay close to the solution of the original problem, if one allows the regularization
parameter $\respar$ to increase slightly.
Apart from the more general, topological setting, the main difference
to~\cite[Lemma~2.2]{Hei08} is the additional inclusion of operator errors into the result.

The next  theorem provides a true stability theorem, including both
data as well as operator perturbations.

\begin{theorem}[Stability]\label{thm:stability:2}
  Let Assumption~\ref{as:stability} hold with $\beta > 0$
  and assume that the inclusion
  \begin{equation}\label{eq:stability2:lim}
    \Phi_\rf(\Fo,y,\respar,t) \subset
      \bigcap_{\delta > 0} \left(\cl \bigcup_{\eps > 0} \Phi_\rf(\Fo,y,\respar-\eps,t+\delta)\right)
  \end{equation}
  holds for every $t \ge 0$.
  Then,
  \begin{equation}\label{eq:stability2:Limsup}
  \emptyset \neq \Limsuptr_{k\to\infty} \sset(\Fo_k,y_k,\respar_k) \subset \sset(\Fo,y,\respar)\;.
  \end{equation}
  If, additionally, the set $\sset(\Fo,y,\respar)$ consists of a single element $x_\respar$,
  then
  \[
    \{x_\respar\} = \Limtr_{k\to\infty} \sset(\Fo_k,y_k,\respar_k)\;.
  \]
\end{theorem}

\begin{proof}
  The convergence of $(\respar_k)_{k\in\N}$ to $\respar$
  and Lemma~\ref{le:Phiincl2} imply that for every $\eps > 0$ and $t \in \R$
  there exists $k_0 \in \N$ such that
  \[
  \Phi_\rf(\Fo,y,\respar-\eps,t)
  \subset \Phi_\rf(\Fo_k,y_k,\respar_k,t)
  \]
  for all $k \ge k_0$.
  Consequently,
  \begin{multline*}
    \limsup_{k\to\infty} v(\Fo_k,y_k,\respar_k)
    = \limsup_{k\to\infty}  \, \Bigl( \inf\set{t:\Phi_\rf(\Fo_k,y_k,\respar_k,t) \neq \emptyset} \Bigr) \\
    \le \inf_{\eps > 0} \, \Bigl(  \inf \set{t:\Phi_\rf(\Fo,y,\respar-\eps,t) \neq\emptyset} \Bigr) \;.
  \end{multline*}
  From~\req{stability2:lim} we obtain that
  \begin{multline*}
  \inf_{\eps > 0} \, \Bigl( \inf \set{t:\Phi_\rf(\Fo,y,\respar-\eps,t) \neq\emptyset} \Bigr)\\
  \le \inf\set{t:\Phi_\rf(\Fo,y,\respar,t)\neq\emptyset} = v(\Fo,y,\respar)\;.
  \end{multline*}
  This shows~\req{stability:V}.
  Now~\req{stability2:Limsup} follows from Lemma~\ref{le:stability}.
\end{proof}

For Theorem~\ref{thm:stability:2} to hold, the mapping
$x \mapsto \sif\bigl(\Fo(x),y\bigr)$ has to satisfy the additional
regularity property~\req{stability2:lim}.
This property requires that every $x \in \X$ for which $\Fo(x) \neq y$
can be approximated by elements $\tilde{x}$ with
$\sif\bigl(\Fo(\tilde{x}),y\bigr) < \sif\bigl(\Fo(x),y\bigr)$
and $\rf(\tilde{x}) \le \rf(x) + \respar$.
That is, the function $x \mapsto \sif\bigl(\Fo(x),y\bigr)$
does not have local minima in the sets $\set{x \in \X:\rf(x)< t}$.
As will be shown in the following Section~\ref{se:linear},
this property is naturally satisfied for linear operators on
Banach spaces.

\subsection{Convergence}

The following theorem  states  the solutions obtained with the
residual method indeed converge to the $\rf$-minimizing solution
of the equation $\Fo(x) = y$,  if  the noise level decreases to zero.
Recall that the set of all   $\rf$-minimizing solution of the equation $F(x) = y$
is given by   $\sset(\Fo,y,0)$.

\begin{theorem}[Convergence]\label{thm:convergence}
  Let  $y \in \Y$ be such that there exists  $x \in \X$ with   $\Fo(x) = y$
  and $\rf(x) < \infty$ and assume that  $\Phi_\rf(\Fo,w,\gamma, t)$ is
  $\tau$-compact for all  $w \in \Y$ and  $\gamma, t \geq 0$.
  If $(y_k)_{k\in\N}$ converges  $\sif$-uniformly to   $y$ and satisfies $\sif( y, y_k) \le \respar_k  \to 0$, then
  \begin{equation}\label{eq:convergence_V}
  \limsup_{k\to\infty} v(\Fo,y_k,\respar_k) \le v(\Fo,y,0) < \infty\;.
  \end{equation}
  In particular,
  \begin{equation}\label{eq:convergence}
  \emptyset \neq \Limsuptr_{k\to\infty} \sset(\Fo,y_k,\respar_k) \subset \sset(\Fo,y,0)\;.
  \end{equation}
  If, additionally, the $\rf$-minimizing solution $x^\dagger$ is unique, then
  \begin{equation}\label{eq:convergence_2}
    \{x^\dagger\} = \Limsuptr_{k\to\infty} \sset(\Fo,y_k,\respar_k)\;.
  \end{equation}
\end{theorem}

\begin{proof}
  By assumption $S(y,y_k) \le \respar_k$,
  which implies that $v(\Fo,y_k,\respar_k) \le \rf(x')$ for all
  $x'\in \Phi(\Fo,y,0)$.
  This proves~\req{convergence_V}.
  Now~\req{convergence} and~\req{convergence_2} follow from Lemma~\ref{le:stability}.
\end{proof}

\section{Linear Spaces}\label{se:linear}

Now we assume that $\X$ and $\Y$ are subsets of
topological vector spaces.
Then the linear structures allows us to introduce more tangible
conditions implying stability of the residual method.

For the following we assume that  $\Fo\colon \X \to \Y$ and $y \in \Y$
are fixed.

\begin{assumption}\label{as:linear}
  Assume that the following hold:
  \begin{enumerate}
  \item\label{it:linear:X} The set $\X$ is a convex subset of a topological
  vector space,  and $\Y$ is a topological vector space.

  \item\label{it:linear:S}
    The mapping $x \mapsto \sif\bigl(\Fo(x_0),y\bigr)$ is semi-strictly quasi-convex.
    That is,
    for all $x_0$,~$x_1 \in \X$ with $\sif\bigl(\Fo(x_0),y\bigr)$, $\sif\bigl(\Fo(x_1),y\bigr) < \infty$,
    and all $0 < \lambda < 1$ we have
    \[
    \sif\bigl(\Fo(\lambda x_0 + (1-\lambda) x_1),y\bigr)
    \le \max\bigl\{ \sif\bigl(\Fo(x_0),y\bigr), \sif\bigl(\Fo(x_1),y\bigr)\bigr\}\;.
    \]
    Moreover, the inequality is strict whenever $\sif\bigl(\Fo(x_0),y\bigr) \neq \sif\bigl(\Fo(x_1),y\bigr)$.

  \item\label{it:linear:domain} For every $\respar \ge 0$ there exists $x \in \X$
    with $\sif\bigl(\Fo(x),y\bigr) \le \respar$ and $\rf(x) < \infty$.

  \item\label{it:linear:R} The domain
    $\domain\rf = \set{x \in \X:\rf(x) < +\infty}$ of $\rf$ is convex and
    for every $x_0,\,x_1 \in \domain\rf$
    the restriction of $\rf$ to
    \[
    L = \set{ \lambda x_0 + (1-\lambda) x_1:0 \le \lambda \le 1}
    \] is continuous.
  \end{enumerate}
\end{assumption}

We now show that Assumption~\ref{as:linear}~implies the main condition
of the stability result Theorem~\ref{thm:stability:2},
the inclusion~\req{stability2:lim}:

\begin{lemma}\label{le:linear:1}
  Assume that Assumption~\ref{as:linear} holds.
  Then~\req{stability2:lim} is satisfied.
\end{lemma}

\begin{proof}
  Let $x_0 \in \Phi_\rf(\Fo,y,\respar,t)$ for some $\respar > 0$.
  We have to show that for every neighborhood $N \subset \X$ of $x_0$ and every
  $\delta > 0$   there exist $\eps > 0$ and $x' \in N$ such that
  $x' \in \Phi_\rf(\Fo,y,\respar-\eps,t+\delta)$.

  Item~\ref{it:linear:domain} in Assumption~\ref{as:linear} implies the
  existence of some $x_1 \in \X$ satisfying the inequalities
  $\sif\bigl(\Fo(x_1),y\bigr) \le \respar/2$ and $\rf(x_1) < \infty$.
  Since we have $\sif\bigl(\Fo(x_1),y\bigr) < \respar$ and $\sif\bigl(\Fo(x_0),y\bigr) \le \respar$,
  we obtain from Item~\ref{it:linear:S}
  that $\sif(\Fo(x),y) < \respar$ for every
  $x \in L := \set{\lambda x_0 + (1-\lambda)x_1:0 \le \lambda < 1}$.
  Since $x_0$,~$x_1 \in \domain\rf$,
  it follows from Item~\ref{it:linear:R} that $\rf$ is continuous on $L$.
  Consequently $\lim_{\lambda \to 1} \rf(\lambda x_0 + (1-\lambda) x_1) = \rf(x_0) \le t$.
  In particular, there exists $\lambda_0 < 1$ such that
  $\rf(\lambda x_0 + (1-\lambda)x_1) \le t+\delta$ for all $1 > \lambda > \lambda_0$.
  Since $\X$ is a topological vector space (Item~\ref{it:linear:X}), it follows that
  $x':=\lambda x_0 + (1-\lambda) x_1 \in N$ for some $1 > \lambda > \lambda_0$.
  This shows the assertion with $\eps := \respar - \sif\bigl(\Fo(x'),y\bigr) > 0$.
\end{proof}

Lemma~\ref{le:linear:1} allows us to apply the stability result Theorem~\ref{thm:stability:2},
which shows that Assumption~\ref{as:linear} implies the continuous dependence
of the solutions of~\req{min_prob} on the data $y$ and the regularization parameter $\respar$.

\begin{proposition}[Stability \& Convergence]\label{pr:linear_stab}
  Let Assumption~\ref{as:linear} hold and
  assume that the sets
  $\Phi_\rf(\Fo,w,\gamma,t)$ are compact for every $\gamma \in \R$,
  $t \in \R$, and $w \in \Y$.
  Assume moreover that $(y_k)_{k\in\N}$ converges $\sif$-uniformly to $y \in \Y$,
  and that $\respar_k \to \respar$.
  If $\respar = 0$, assume in addition that $\sif(y,y_k) \le \respar_k$.
  Then
  \[
  \emptyset \neq \Limsuptr_{k\to\infty} \sset(\Fo,y_k,\respar_k) \subset \sset(\Fo,y,\respar)\;.
  \]
  If, additionally, the set $\sset(\Fo,y,\respar)$ consists of a single element $x_\respar$,
  then
  \[
  \{x_\respar\} = \Limsuptr_{k\to\infty}  \sset(\Fo,y_k,\respar_k)\;.
  \]
\end{proposition}

\begin{proof}
  If $\respar = 0$, the assertion follows from Theorem~\ref{thm:convergence}.
  In the case $\respar > 0$, Lemma~\ref{le:linear:1}
  implies that~\req{stability2:lim} holds.
  Thus, the assertion follows from Theorem~\ref{thm:stability:2}.
  Note that the non-emptyness of the sets $\Phi_\rf(\Fo,w,\gamma,t)$
  for some $t$ follows from Item \ref{it:linear:domain}
  in Assumption \ref{as:linear}.
\end{proof}

\begin{proposition}[Stability]\label{pr:linear_stab2}
  Let Assumption~\ref{as:linear} hold.
  Assume that $(y_k)_{k\in\N}$ converges $\sif$-uniformly to $y \in \Y$,
  the mappings $\Fo_k\colon \X \to \Y$ converge locally $\sif$-uniformly
  to $\Fo\colon \X \to \Y$ (see Definition~\ref{de:s_uniform}), and $\respar_k \to \respar > 0$.
  Assume that the sets $\Phi_\rf(\Fo_k,w,\gamma,t)$ and $\Phi_\rf(\Fo,w,\gamma,t)$ are compact for every
  $\gamma \ge 0$, $t \in \R$, and $w \in \Y$.
  Then
  \[
  \emptyset \neq \Limsuptr_{k\to\infty} \sset(\Fo_k,y_k,\respar_k) \subset \sset(\Fo,y,\respar)\;.
  \]
  If, additionally, the set $\sset(\Fo,y,\respar)$ consists of a single element $x_\respar$,
  then
  \[
    \{x_\respar\} = \Limtr_{k\to\infty} \sset(\Fo_k,y_k,\respar_k)\;.
  \]
\end{proposition}

\begin{proof}
  Again, Lemma~\ref{le:linear:1} shows that~\req{stability2:lim} holds.
  Moreover, the non-empty\-ness of the sets $\Phi_\rf(\Fo,w,\gamma,t)$
  and $\Phi_\rf(\Fo_k,w,\gamma,t)$ (at least for $k$ sufficiently large)
  for some $t$ follows from Item \ref{it:linear:domain}
  in Assumption \ref{as:linear} and the local $\sif$-uniform convergence
  of the mappings $\Fo_k$ to $\Fo$.
  Thus the assertion follows from Theorem~\ref{thm:stability:2}.
\end{proof}

Item~\ref{it:linear:S} in Assumption~\ref{as:linear} is concerned with the interplay
of the functional $\Fo$ and the distance measure $\sif$.
The next two examples consider two situations, where this part of the assumption holds.
Example~\ref{ex:linear1} considers linear operators $\Fo$ and convex distance measures $\sif$.
Example~\ref{ex:stab_nonlinear} introduces a class of \emph{non-linear} operators on Hilbert spaces,
where Item~\ref{it:linear:S} is satisfied if the distance measure equals
the squared Hilbert space norm.

\begin{example}\label{ex:linear1}
  Assume that $\Fo \colon \X \to \Y$ is linear and $\sif$
  is convex in its first component.
  Then Item~\ref{it:linear:S} in Assumption~\ref{as:linear} is satisfied.
  Indeed, in such a situation,
  \begin{align*}
  \sif\bigl(\Fo(\lambda x_0 + (1-\lambda)x_1),y\bigr)
  &= \sif\bigl(\lambda \Fo(x_0) + (1-\lambda)\Fo(x_1),y\bigr)
  \\
  &\le \lambda \sif\bigl(\Fo(x_0),y\bigr) + (1-\lambda)\sif\bigl(\Fo(x_1),y\bigr)
  \\&\le
  \max\bigl\{\sif\bigl(\Fo(x_0),y\bigr), \sif\bigr(\Fo(x_1),y\bigr)\bigr\}\;.
  \end{align*}
  If moreover, $\sif\bigl(\Fo(x_0),y\bigr) \neq \sif\bigr(\Fo(x_1),y\bigr)$
  and $0 < \lambda < 1$, then the last inequality is strict.
\end{example}

\begin{example}\label{ex:stab_nonlinear}
  Assume that $\Y$ is a Hilbert space, $\sif(y,z) = \norm{y-z}^2$,
  and $\Fo \colon \X \to \Y$ is two times G\^ateaux differentiable.
  Then Item~\ref{it:linear:S} in Assumption~\ref{as:linear}
  holds if for all $x_0 \neq x_1 \in \X$ the mapping
  \[
  t \mapsto \varphi(t;x_0,x_1) := \norm{\Fo(x_0 + tx_1) - y}^2
  \]
  has no local maxima.
  This condition holds, if the inequality $\partial_t^2 \varphi(0;x_0,x_1) > 0$ is
  satisfied whenever $\partial_t \varphi(0;x_0,x_1) = 0$.
  The computation of the derivative  of $\varphi(\,\cdot\,;x_0,x_1)$ at zero yields that
  \[
  \partial_t \varphi(0;x_0,x_1) = 2\inner{\Fo'(x_0)(x_1)}{\Fo(x_0)}
  \]
  and
  \[
  \partial_t^2 \varphi(0;x_0,x_1) = 2\inner{\Fo''(x_0) (x_1;x_1)}{\Fo(x_0)} + 2\norm{\Fo'(x_0)x_1}^2\;.
  \]

  Consequently, Item~\ref{it:linear:S} in Assumption~\ref{as:linear} is satisfied
  if, for every $x_0$, $x_1 \in \X$ with $x_1 \neq 0$, the equality
  $\inner{\Fo'(x_0)(x_1)}{\Fo(x_0)} = 0$ implies that
  \[
    \inner{\Fo''(x_0) (x_1;x_1)}{\Fo(x_0)} + \norm{\Fo'(x_0)(x_1)}^2 > 0\;.
  \]
\end{example}

\subsection{Regularization on $L^p$-spaces}

Let $p \in (1,\infty)$ and set $\X = L^p(\Omega,\mu)$ for some $\sigma$-finite measure space
$(\Omega,\mu)$.
  Assume that $\Y$ is a Banach space and $\Fo\colon \X \to \Y$ is
  a bounded linear operator with dense range.
  Let $\rf(x) = \norm{x}_p^p$ and $\sif(w,y) = \norm{w-y}$.
  We thus consider the minimization problem
  \[
  \norm{x}_p^p \to \min
  \qquad\text{ subject to }\quad
  \norm{\Fo x - y} \le \respar\;.
  \]

  We now show that in this situation the assumptions of Proposition~\ref{pr:linear_stab}
  are satisfied.
  To that end, let $\tau$ be the weak topology on $L^p(\Omega,\mu)$.
  As $L^p(\Omega,\mu)$ is reflexive, the level sets $\set{x \in \X:\rf(x) \le t}$ are weakly compact.
  Moreover, the mapping $x \mapsto \norm{\Fo x-y}$ is weakly lower semi-continuous.
  Thus all the sets $\Phi_\rf(\Fo,w,\gamma,t)$ are weakly compact.
  Example~\ref{ex:linear1} shows that Item~\ref{it:linear:S} in Assumption~\ref{as:linear} holds.
  Item~\ref{it:linear:domain} follows from the density of the range of $\Fo$.
  Finally, Item~\ref{it:linear:R} holds, because $\rf$ is norm continuous
  and convex.

  Now assume that $y_k \to y$ and $\respar_k \to \respar$.
  If $\respar = 0$ assume in addition that $\norm{y_k-y} \le \respar_k$.
  The strict convexity of $\rf$ and convexity of the mappings $x \mapsto \norm{\Fo x-y_k}$
  imply that each set $\sset(\Fo,y_k,\respar_k)$ consists of a single element $x_k$.
  Similarly, $\sset(\Fo,y,\respar)$ consists of a single element $x^\dagger$.
  From Proposition~\ref{pr:linear_stab} we now obtain that $(x_k)_{k\in\N}$ weakly converges to $x^\dagger$
  and $\norm{x_k}_p^p \to \norm{x^\dagger}_p^p$.
  Thus, in fact, the sequence $(x_k)_{k\in\N}$ strongly converges to $x^\dagger$
  (see~\cite[Cor.~5.2.19]{Meg98}).

  Let $\respar > 0$ and assume that $\Fo_k\colon \X \to \Y$
  is a sequence of bounded linear operators converging to $\Fo$
  with respect to the strong topology on $L(\X,\Y)$,
  that is, $\sup\{\norm{\Fo_k x - \Fo x}:\norm{x} \le 1\} \to 0$.
  Let again $\respar_k \to \respar$ and $y_k \to y$, and denote
  by $x_k$ the single element in $\sset(\Fo_k,y_k,\respar_k)$.
  Applying Proposition~\ref{pr:linear_stab2},
  we again obtain that $x_k \to x^\dagger$.
%\end{example}

\begin{remark}
  The above results rely heavily on the assumption that $p > 1$,
  which implies that the space $L^p(\Omega,\mu)$ is reflexive.
  In the case $\X = L^1(\Omega,\mu)$, the level sets $\set{x \in \X:\norm{x}_1 \le t}$
  fail to be weakly compact, and thus even the existence of a solution
  of Problem~\req{min_prob} need not hold.
\end{remark}

\begin{remark}
  The assertions concerning stability and convergence
  with respect to the norm topology remain valid, if $\X$ is any uniformly
  convex Banach space and $\rf$ the norm on $\X$ to some power $p > 1$.
  Also in this case, weak convergence and convergence of norms imply
  the strong convergence of a sequence \cite[Thm.~5.2.18]{Meg98}.
  More generally, this property is called the \emph{Radon--Riesz property} \cite[p.~453]{Meg98}.
  Spaces satisfying this property are also called \emph{Efimov--Stechkin} spaces in~\cite{Tan97}.
\end{remark}

\subsection{Regularization of Probability Measures}
\label{ss:prob}

Let $(\Omega, d)$ be a separable, complete metric space with distance $d$ and denote by
$\mathcal{P}(\Omega)$ the space of probability measures on the Borel sets of $\Omega$.
That is, $\mathcal{P}(\Omega)$ consists of all positive Borel measures $\mu$
on $\Omega$ that satisfy $\mu(\Omega) = 1$.
For $p \ge 1$ the \emph{$p$-Wasserstein distance}
on $\mathcal{P}(\Omega)$ is defined as
\[
W_p(\mu,\nu) := \left(\inf\set{\int d(x,y)^p\,d\xi:\xi \in \mathcal{P}(\Omega\times\Omega),\
\pi^1_\# \xi = \mu,\ \pi^2_\# \xi = \nu}\right)^{1/p}.
\]
Here $\pi^i_\# \xi$ denotes the push forward of the measure $\xi$ by means of the $i$-th projection.
In other words, $\pi^1_\# \xi (U) = \xi(U\times \Omega)$
and $\pi^2_\#\xi(U) = \xi(\Omega\times U)$ for every Borel set $U \subset \Omega$.

Recall that the \emph{narrow topology} on $\mathcal{P}(\Omega)$ is induced
by the action of elements of $\mathcal{P}(\Omega)$ on continuous functions $u \in C(\Omega)$.
That is, a sequence $(\mu_k)_{k\in\N}\subset \mathcal{P}(\Omega)$ converges
narrowly to $\mu\in\mathcal{P}(\Omega)$, if
\[
    \int_\Omega u\,d\mu_k = \int_\Omega u\,d\mu  \quad \text{ for all } u \in C(\Omega)  \,.
\]

\begin{lemma}\label{le:wp}
  Let $p \ge 1$.
  Then the Wasserstein distance satisfies, for every $\mu_1$, $\mu_2$, $\nu \in \mathcal{P}(\Omega)$
  and $0 \le \lambda \le 1$, the inequality
  \begin{equation}\label{eq:wp}
  W_p\bigl(\lambda\mu_1 + (1-\lambda)\mu_2,\nu\bigr)^p
  \le \lambda W_p(\mu_1,\nu)^p + (1-\lambda)W_p(\mu_2,\nu)^p\;.
  \end{equation}
  Moreover it is lower semi-continuous with respect to the narrow topology.
\end{lemma}

\begin{proof}
  The lower semi-continuity of $W_p$ has, for instance, been shown in~\cite{GivSho84}.
  In order to show the inequality~\req{wp},
  let $\xi_1$, $\xi_2 \in \mathcal{P}(\Omega\times\Omega)$ be
  two measures that realize the infimum in the definition
  of $W_p(\mu_1,\nu)$ and $W_p(\mu_2,\nu)$, respectively.
  Then $\pi^1_\#\bigl(\lambda\xi_1+(1-\lambda)\xi_2\bigr) = \lambda\mu_1 + (1-\lambda)\mu_2$
  and $\pi^2_\#\bigl(\lambda\xi_1+(1-\lambda)\xi_2\bigr) = \nu$,
  which implies that the measure $\lambda\xi_1+(1-\lambda)\xi_2$ is admissible
  for measuring the distance between $\lambda\mu_1+(1-\lambda)\mu_2$ and $\nu$.
  Therefore
  \[
  \begin{aligned}
    W_p\bigl(\lambda\mu_1 + & (1-\lambda)\mu_2,\nu\bigr)^p
     \\
     &= \inf\set{\int d(x,y)^p\,d\xi:\pi^1_\# \xi = \lambda\mu_1+(1-\lambda)\mu_2,\ \pi^2_\# \xi = \nu}\\
    & \le \int d(x,y)^p\,d(\lambda\xi_1 + (1-\lambda)\xi_2)\\
    & = \lambda W_p(\mu_1,\nu)^p + (1-\lambda) W_p(\mu_2,\nu)^p\,,
  \end{aligned}
  \]
  which proves the assertion.
\end{proof}

Since $\mathcal{P}(\Omega)$ is a convex subset of the space $\mathcal{M}(\Omega)$
of all finite Radon measures on $\Omega$,
and the narrow topology on $\mathcal{P}(\Omega)$ is the restriction
of the weak$^*$ topology on $\mathcal{M}(\Omega)$ considered as the dual of
$C_b(\Omega)$, the space of bounded continuous functions on $\Omega$,
it is possible to apply the results of this section also to the situation
where $\Y = \mathcal{P}(\Omega)$ and $\sif = W_p$.
As an easy example, we consider the problem of density estimation from
a finite number of measurements.

\begin{example}
  Let $\Omega \subset \R^n$ be an open domain.
  Given a finite number of measurements $\{y_1,\ldots,y_k\}\subset \Omega$,
  the task of density estimation is the problem of finding a simple density
  function $u$ on $\Omega$ in such a way that the measurements look like a typical sample
  of the distribution defined by $u$.
  Interpreting the measurements as a normalized sum of delta peaks,
  that is, equating $\{y_1,\ldots,y_k\}$ with the measure
  $\mathbf{y} := \frac{1}{k}\sum_i \delta(y_i) \in \mathcal{P}(\Omega)$,
  we can easily translate the problem into the setting of this paper.

  We set $\X := \set{u\in L^1(\Omega):u \ge 0 \text{ and } \norm{u}_1 = 1}$,
  which is a convex and closed subset of $L^1(\Omega)$,
  $\Y := \mathcal{P}(\Omega)$, and consider the embedding $\Fo\colon \X\to \mathcal{P}(\Omega)$,
  $u \mapsto u\,\mathcal{L}^n$.
  Then $\Fo$ is continuous with respect to the weak topology on $\X$
  and the narrow topology on $\mathcal{P}(\Omega)$.
  We now consider the distance measure $\sif = W_p$ for some $p \ge 1$
  and the Euclidean distance $d$ on $\Omega$.
  Then Lemma~\ref{le:wp} implies that, for every $\mu \in \mathcal{P}(\Omega)$,
  the mapping $u \mapsto W_p(Fu,\mu)$ is weakly lower semi-continuous.

  There are several possibilities for choosing a regularization functional on $\X$.
  If $\Omega$ is bounded (or at least $\mathcal{L}^n(\Omega) < \infty$),
  one can, for instance, use the Boltzmann--Shannon entropy defined by
  \[
  \rf(u) := \int_\Omega u\,\log(u)\,dx
  \qquad
  \text{ for } u \in \X\;.
  \]
  Then the Theorems of De la Vall\'ee Poussin and Dunford--Pettis
  (see~\cite[Thms.~2.29, 2.54]{FonLeo07}) show
  that the lower level sets of $\rf$ are weakly pre-compact in $L^1(\Omega)$.
  Moreover, the functional $\rf$ is convex and therefore weakly lower semi-continuous
  (see~\cite[Thm.~5.14]{FonLeo07}).
  Using Proposition~\ref{prop:compact}, we therefore obtain that
  the compactness required in Assumption~\ref{as:stability} holds.
  Also, Lemma~\ref{le:wp} shows that Item~\ref{it:linear:S} in Assumption~\ref{as:linear} holds.
  Items~\ref{it:linear:X} and~\ref{it:linear:R} are trivially satisfied.
  Finally, Item~\ref{it:linear:domain} follows from the density of $\domain{R}$ in $\X$
  and the density (with respect to the narrow topology) of $\range \Fo$ in $\mathcal{P}(\Omega)$.
  In addition, it has been shown in~\cite{Vis84} that the weak convergence
  of a sequence $(u_k)_{k\in\N} \subset L^1(\Omega)$ to $u \in L^1(\Omega)$
  together with the convergence $\rf(u_k) \to \rf(u)$
  imply that $\norm{u_k-u}_1 \to 0$.
  Thus the topology $\tau_{\rf}$ coincides with the
  strong topology on $\X$.

  Proposition~\ref{pr:linear_stab} therefore implies that the residual method is
  a stable and convergent regularization method with respect to the strong topology on $\X$.
  More precisely, given a sample $\mathbf{y} = \frac{1}{k} \sum_i \delta(y_i)$, the density estimate $u$
  depends continuously on the positions $y_i$ of the measurements and on the regularization
  parameter $\respar$.
  In addition, if the number of measurements increases, then
  the Wasserstein distance between the sample and the true probability
  converges almost surely to zero.
  Thus also the reconstructed density converges to the true underlying
  density, provided the regularization parameters decrease to zero slowly enough.
\end{example}

\section{Convergence Rates}\label{se:rates}

In this section we derive quantitative estimates (convergence rates)
for the difference between regularized solutions $x_\respar \in \sset(\Fo,y,\respar)$
and the exact solution of the equation $\Fo(x^\dagger) = y^\dagger$.

\medskip
For Tikhonov regularization, convergence rates have been derived
in~\cite{BotHof10,BurOsh04,FleHof10,HofKalPoeSch07,Res05,ResSch06}
in terms of the \emph{Bregman distance}. However, its classical definition,
  \begin{equation}\label{eq:bregman_dist}
  D_\xi(x,x^\dagger) = \rf(x) - \rf(x^\dagger)  + \minner{\xi}{x^\dagger - x}_{\X^*,\X} \,,\end{equation}
where $\xi \in \partial\rf(x^\dagger) \subset \X^*$,
requires the space $\X$ to be linear and the functional $\rf$ to be convex, as the
(standard) subdifferential $\partial\rf(x^\dagger)$ is only defined for convex functionals.
In the sequel we will extend the notion of subdifferentials and Bregman distances to work
for arbitrary functionals $\rf$ on arbitrary sets $\X$.
To that end, we make use of a generalized notion of convexity, which is not
based on the duality between a Banach space $\X$ and its dual $\X^*$
but on more general pairings (see~\cite{Sin97}).
The same notion has recently been used in~\cite{Gra10b} for the derivation of
convergence rates for non-convex regularization functionals.

\begin{definition}[Generalized Bregman Distance]\label{def:bregman}
  Let $W$ be a set of functions $w\colon \X \to \R$, let
  $\rf\colon \X \to \R\cup\{+\infty\}$ be  a functional and let
  $x^\dag \in X$.

  \begin{enumerate}
  \item[(a)]
   The  functional $\rf$ is \emph{convex at $x^\dag$ with respect to $W$}, if
   \begin{equation}\label{eq:conj}
   \rf(x^\dag ) = \rf^{**}(x^\dag)
   :=  \sup_{w \in W} \,
   \Bigl( \inf_{x\in \X}  \bigl( \rf(x) - w(x) + w(x^\dag) \bigr) \Bigr) \,.
   \end{equation}

  \item[(b)]
 Let  $\rf$ be convex at $x^\dag$ with respect to $W$.
  The \emph{subdifferential at $x^\dag$ with respect $W$}  is defined as
  \[
  \partial_W \rf(x^\dag)
  := \set{w \in W: \rf(x) \ge \rf(x^\dag) + w \bigl(x\bigr) - w(x^\dag)
    \text{ for all } x \in \X}.
  \]

  \item[(c)]
  Let  $\rf$ be convex at $x^\dag$ with respect to $W$.
  For $w \in \partial_W \rf(x^\dag)$ and $x \in X$, the
  \emph{Bregman distance between $x^\dag$ and $x$ with respect to $w$}
  is defined as
  \begin{equation}\label{eq:bregman_dist2}
    D_w (x,x^\dag) := \rf(x) - \rf(x^\dag) - w(x) + w(x^\dag) \,.
  \end{equation}
  \end{enumerate}
\end{definition}

\begin{remark}\label{rem:bregman-lc}
  Let $\X$ be a Banach space and set $W = \X^*$.
  Then a functional $\rf \colon  \X  \to  \R\cup\{+\infty\}$
  is convex with respect to $W$,
  if and only if it is lower semi-continuous and convex in the classical sense.
  Moreover,at  every $x^\dag \in \X$, the subdifferential with respect $W$ coincides with
  the classical subdifferential  $\partial\rf(x^\dagger) \subset \X^*$.
  Finally, the standard Bregman distance, defined by \req{bregman_dist},
  coincides with the Bregman distance obtained by means of Definition~\ref{def:bregman}.
\end{remark}

In the following, let $W$ be a given family of real valued functions on
$\X$. Convergence rates in Bregman distance with respect to $W$
will be derived under the following assumption:

\begin{assumption}\label{as:rates} \mbox{}
\begin{enumerate}
\item
  There exists a monotonically increasing function
  $\psi \colon [0, \infty) \to [0, \infty) $ such that
  \begin{equation}\label{eq:triangle}
        \sif(y_1,y_2) \le \psi \left( \sif(y_1,y_3) + \sif(y_2,y_3) \right)
        \qquad  \text{ for all } y_1, y_2, y_3 \in \Y  \,.
  \end{equation}

\item
For some given point   $x^\dag \in \X$, the functional 
$\rf\colon \X \to \R\cup\{+\infty\}$ is convex at $x^\dag$ 
with respect to $W$.

\item
There exist  $w \in \partial_W \rf(x^\dagger)$
and constants $\gamma_1 \in [0,1)$, $\gamma_2 \ge 0$ such that
  \begin{equation}\label{eq:bregman_est}
    w(x^\dagger)-w(x)
    \le \gamma_1 D_w(x,x^\dagger) + \gamma_2 \sif\bigl( \Fo(x), \Fo(x^\dagger) \bigr)
  \end{equation}
    for every $x \in \Phi_\rf\left(\Fo,\Fo(x^\dagger),\psi(2\respar),\rf(x^\dagger)\right)$.
  \end{enumerate}
\end{assumption}

In a Banach space setting, the \emph{source  inequality}~\req{bregman_est} has already been used
in~\cite{HofKalPoeSch07,SchGraGroHalLen09} to derive convergence rates for Tikhonov regularization with convex
functionals and in \cite{Hei08} for multiparameter regularization.
Equation \req{triangle} is an alternate for the missing triangle inequality
in the non-metric case.

\begin{theorem}[Convergence Rates]\label{thm:rates}
  Let Assumption~\ref{as:rates} hold and let $y \in \Y$ satisfy
  $\sif \left(\Fo(x^\dagger),y\right) \le \respar$.
  Then, the   estimate
  \begin{equation}\label{eq:rate_est}
    D_w(x_\respar, x^\dagger)
    \le
    \frac{\gamma_2 }{1-\gamma_1} \, \psi\bigl(\respar + \sif\bigl(\Fo(x^\dagger),y\bigr)\bigr)
  \end{equation}
   holds for all  $x_\respar \in \sset(\Fo,y,\respar)$.
\end{theorem}

\begin{proof}
  Let $x_\respar \in \sset(\Fo,y,\respar)$.
  This, together with~\req{triangle} and the assumption that
  $\sif\left(\Fo(x^\dagger),y\right)\le\respar$,
  implies that
  \[
  \sif\bigl(\Fo(x_\respar),\Fo(x^\dagger)\bigr)
  \le \psi  \bigl( \sif\bigl(\Fo(x_\respar),y\bigr) + \sif\bigl(\Fo(x^\dag),y\bigr) \bigr)
  \le \psi(2\respar)\,.
  \]
  Together with~\req{bregman_est} it follows that
   \begin{multline*}
    D_w(x_\respar,x^\dagger)
    = \rf(x_\respar) - \rf(x^\dagger) - w(x_\respar) + w(x^\dagger)\\
    \le \rf(x_\respar) - \rf(x^\dagger) + \gamma_1 D_w(x_\respar,x^\dagger)
       + \gamma_2 \sif\bigl(\Fo(x_\respar),\Fo(x^\dagger)\bigr)\;.
  \end{multline*}
  The assumption $\gamma_1 \in [0,1)$ implies the  inequality
  \begin{equation} \label{eq:rates:hlp1}
    D_w(x_\respar,x^\dagger)
    \leq \frac{1}{1-\gamma_1}\bigl(\rf(x_\respar) - \rf(x^\dagger)\bigr)
           + \frac{\gamma_2}{1-\gamma_1}\,\sif\bigl(\Fo(x_\respar),\Fo(x^\dagger)\bigr)\;.
  \end{equation}
  Since $\sif(\Fo(x^\dagger),y ) \le \respar$, we have
  $ \rf(x_\respar) \leq \rf(x^\dagger)$. Therefore \req{rates:hlp1} and~\req{triangle} imply
  \begin{equation*}
      D_w(x_\respar, x^\dagger)
    \le \frac{\gamma_2}{1-\gamma_1}\,\sif\bigl(\Fo(x_\respar),\Fo(x^\dagger)\bigr)
    \leq \frac{\gamma_2}{1-\gamma_1} \, \psi\bigl(\respar + \sif(\Fo(x^\dagger),y)\bigr) \;,
  \end{equation*}
  which concludes the proof.
\end{proof}

\begin{remark}
  Typically, convergence rates are formulated in a setting  which slightly differs from
  the one of Theorem \ref{thm:rates},
  see \cite{BurOsh04,EngHanNeu96,HofKalPoeSch07,SchGraGroHalLen09}.
  There one assumes the existence of an $\rf$-minimizing solution $x^\dagger \in \X$
  of the equation $\Fo(x^\dagger) = y^\dagger$, for some exact data
  $y^\dagger \in \range(\Fo)$.
  Instead of $y^\dagger$, only noisy data  $y \in \Y$
  and the error bound $\sif(y^\dagger,y) \le \beta$ are given.

  For this setting, \req{rate_est} implies the rate
  \[
  D_w( x_\respar, x^\dagger)
  \le \frac{\gamma_2}{1-\gamma_1} \, \psi(2 \respar)
  = \mathcal{O}\bigl( \psi(2\respar) \bigr)
  \qquad
  \text{ as } \respar \to 0 \,,
  \]
  where $x_\respar \in \sset(\Fo,y,\respar)$ denotes any regularized solution.
\end{remark}

\begin{remark}\label{re:bregman_est_var}
  The inequality~\req{bregman_est} is equivalent to the existence
  of $\eta_1$, $\eta_2 > 0$ such that
  \begin{equation}\label{eq:bregman_est_var}
  w(x^\dagger) - w(x)
  \le \eta_1 \bigl(\rf(x)-\rf(x^\dagger)\bigr) + \eta_2 \, \sif\bigl(\Fo(x),\Fo(x^\dagger)\bigr)\;.
  \end{equation}
  Indeed, we obtain~\req{bregman_est_var} from~\req{bregman_est}
  by setting $\eta_1 := \gamma_1/(1-\gamma_1)$ and $\eta_2 := \gamma_2/(1-\gamma_1)$.
  Conversely, \req{bregman_est_var} implies~\req{bregman_est} by taking
  $\gamma_1 := \eta_1/(1+\eta_1)$ and $\gamma_2 := \eta_2/(1+\eta_1)$.
\end{remark}

\subsection{Convergence Rates in Banach spaces}

In the following, assume that $\X$ and $\Y$ are Banach spaces with norms
$\enorm$ and $\enorm$, and assume that
 $\rf$ is a   convex and lower semi-continuous functional on $\X$.
We set $\sif(y,z):= \norm{y-z}$  and  let
$D_{\xi}$ with $\xi \in \partial\rf(x^\dagger)$
denote the classical Bregman distance  (see Remark \ref{rem:bregman-lc}).

\medskip
If $x^\dagger$ satisfies the   inequality
\begin{equation}\label{eq:bregman_est_banach}
  \minner{\xi}{x^\dag-x}
  \le \gamma_1 D_{\xi}(x,x^\dagger) +
  \gamma_2 \bigl\lVert  \Fo(x) - \Fo(x^\dagger) \bigr\rVert 
\end{equation}
and  $y$ are given data with  $\snorm{\Fo(x^\dagger) - y} \leq \respar$,
then Theorem~\ref{thm:rates} implies the convergence
rate $D_\xi(x_\respar, x^\dagger) = \mathcal O(\respar)$.
In the special case where $\X$ is a Hilbert space and $\rf(x) =  \norm{x}^2$
we have $D_\xi(x, x^\dagger) =   \snorm{x-x^\dagger}^2$, which implies the
convergence  rate $\snorm{x-x^\dagger} = \mathcal O\bigl(\respar^{1/2}\bigr)$ with respect to  the norm.
In Proposition~\ref{pr:rates_norm} below we show that the same convergence rate
holds on any $2$-convex space.
For $r$-convex Banach spaces with $r > 2$, we derive the rate $\mathcal{O}\bigl(\respar^{1/r}\bigr)$.

\begin{definition}
The  Banach  space $\X$ is called $r$-convex (or is said to have modulus of convexity
of power type $r$),  if there exists a
constant $C>0$ such that
 \[
    \inf \set{1- \norm{(x+y)/2} :\norm{x} = \norm {y} =1,\, \norm{x-y}\ge \epsilon} \ge C \eps^r
\]
for all $\eps \in [0,2]$.
\end{definition}

Note that every  Hilbert space is  $2$-convex and
that there is \emph{no} Banach space (with $\dim(\X)\ge 2$) that is $r$-convex for
some $r<2$ (see \cite[pp.~63ff]{LinTza79}).

\begin{proposition}[Convergence rates in the norm]\label{pr:rates_norm}
  Let $\X$ be an $r$-convex Banach space with $r \ge 2$ and let
  $\rf(x) := \norm{x}^r/r$.
  Assume that there exists $x^\dagger \in \X$, a subgradient $\xi \in \partial \rf(x^\dagger)$,
  and constants $\gamma_1 \in [0,1)$, $\gamma_2 \ge 0$, $\respar_0 >0$ such that \req{bregman_est_banach}
  holds for  every $x \in \Phi_\rf\bigl(\Fo,\Fo(x^\dagger),2\respar_0,\rf(x^\dagger)\bigr)$.

  Then there exists a constant $c>0$ such that
  \begin{equation}\label{eq:rate_est_banach}
  \mnorm{x_\respar - x^\dagger}
    \leq
    c \bigl(\respar  + \mnorm{\Fo(x^\dagger)- y} \bigr)^{1/r}
  \end{equation}
  for all  $\respar \in  [0, \beta_0]$, all $y \in \Y$ with $\snorm{\Fo(x^\dagger)- y } \le \respar$,  and all
  $x_\respar \in \sset(\Fo,y,\respar)$.
\end{proposition}

\begin{proof}
  Let $J_r: \X \to 2^{\X^*}$ denote the duality mapping with respect to the weight
  function $s \mapsto s^{r-1}$.
  In \cite[Equation $(2.17)'$]{XuRoa91} it is shown that there exists a constant $K > 0$ such that
  \begin{equation}\label{eq:xuroach}
    \mnorm{x^\dagger + z}^r \ge
    \mnorm{x^\dagger}^r + r \minner{j_r(x^\dagger)}{z}_{\X^*,\X}
    +
    K  \norm{z}^r
  \end{equation}
  for all $x^\dagger$, $z \in \X$ and $j_r(x^\dagger) \in J_r(x^\dagger)$.
  By Asplund's theorem \cite[Chap.~1, Thm.~4.4]{Cio90},
  the duality mapping $J_r$ equals the subgradient of   $\rf= \enorm^r/r$.
  Therefore, by taking $z =  x-x^\dagger$ and $j_r(x^\dagger)=\xi$, inequality \req{xuroach}  implies
  \begin{equation}\label{eq:bregam-norm}
    D_\xi(x, x^\dagger) \ge \frac{K}{r} \ \mnorm{x-x^\dagger}^r
    \qquad \text{ for all } x^\dagger, x \in \X  \text{ and }  \xi \in \partial \rf(x^\dagger)\,.
  \end{equation}
  Consequently, \req{rate_est_banach}  follows from Theorem~\ref{thm:rates}.
\end{proof}

Exact values for the constant $K$ in \req{bregam-norm}
(and thus for the constant $c$ in \req{rate_est_banach}) can be derived from \cite{XuRoa91}.
Bregman distances satisfying \req{bregam-norm} are called $r$-coercive in \cite{HeiHof09}.
This $r$-coercivity has already been applied in \cite{BonKazMaaSchoSchu08} for
the minimization of Tikhonov functionals in Banach spaces.\smallskip

\begin{example}
The spaces $\X = L^p(\Omega, \mu)$ for $p \in (1,2]$ and some $\sigma$-finite measure space
$(\Omega,\mu)$ are examples of 2-convex Banach spaces
(see~\cite[p.~81, Remarks following Theorem 1.f.1.]{LinTza79}).
Consequently we obtain for these spaces
the convergence rate $\mathcal O\bigl(\respar^{1/2}\bigr)$.
The spaces $\X = L^p(\Omega, \mu)$ for $p > 2$ are only $p$-convex,
leading to the rate $\mathcal O\bigl(\respar^{1/p}\bigr)$ in those spaces.
\end{example}

\begin{remark}\label{rem:rates:banach}
  The book \cite[pp.~70ff]{SchGraGroHalLen09} clarifies the relation between~\req{bregman_est_banach}
  and the source conditions used to derive convergence rates for convex functionals
  on Banach spaces.
  In particular, it is shown that, if $\Fo$ and $\rf$ are G\^ateaux differentiable at $x^\dagger$
  and there exist $\gamma > 0$ and $\omega \in \Y^*$ such that $\gamma \norm{\omega} < 1$
  and
  \begin{gather}
    \xi = \Fo'(x^\dagger)^*\,\omega \in \partial \mathcal R(x^\dagger) \,, \label{eq:cond-a}
    \\ \label{eq:cond-b}
    \mnorm{\Fo(x) - \Fo(x^\dagger)- \Fo'(x^\dagger)(x-x^\dagger)}  \leq \gamma D_\xi(x,x^\dagger)
  \end{gather}
  for every $x \in \X$, then~\req{bregman_est_banach} holds on $\X$.
  (Here $\Fo'(x^\dagger)^* : \Y^* \to \X^*$ is the adjoint of $\Fo'(x^\dagger)$.)
  Conversely, if $\xi \in \partial \rf(x^\dagger)$ satisfies \req{bregman_est_banach},
  then~\req{cond-a} holds for every $x \in \X$.

  In the particular case that $\Fo\colon \X \to \Y$ is linear and bounded,
  the inequality~\req{cond-b} is trivially satisfied with $\gamma =0$.
  Thus, \req{bregman_est_banach} is equivalent to
  the sourcewise representability of the subgradient,
  $ \xi \in  \partial \rf(x^\dagger) \cap \range(\Fo^*)$.
\end{remark}

\section{Sparse Regularization}
\label{sec:sparsity}

Let $\Lambda$ be an at most countable index set, define
\[
\ell^2(\Lambda) := \biggl\{x = (x_\lambda)_{\lambda \in \Lambda} \subset \R: \sum_{\lambda \in \Lambda} \abs{x_\lambda}^2 < \infty \biggr\}\,,
\]
and assume that $\Fo\colon \X:=\ell^2(\Lambda) \to \Y$ is a bounded linear operator with dense range
in the Hilbert space $\Y$.
We consider for $p \in (0,2)$ the minimization problem
\begin{equation}\label{eq:min_prob_sparse}
  \rf_{p}(x):= \norm{x}_{\ell^p(\Lambda)}^p  := \sum_{\lambda \in \Lambda} \abs{x_\lambda}^p \to \min
  \qquad\text{ subject to }\quad \norm{\Fo x  - y}^2  \le \respar\;.
\end{equation}
For $p>1$, the subdifferential $\partial\rf_p(x^\dagger)$ is at most single valued
and is identified with its single element. 
(The  subdifferential may be empty since we consider   
$\rf_{p}$ as functions on $\ell^2(\Lambda) $.)

\begin{remark}[Compressed Sensing]
In a finite dimensional setting with $p = 1$, the minimization problem
\req{min_prob_sparse} has received a lot of attention during the last years
under the name of \emph{compressed sensing}
(see \cite{Can08,CanRomTao06,CanTao06,Don06c, Don06b,DonElaTem06,Ela10,Fuc05b,Tro06}).
Under some assumptions, the solution of~\req{min_prob_sparse} with $y = \Fo x^\dagger$ and $\respar = 0$
has been shown to recover $x^\dagger$ \emph{exactly}
provided the set
$\{\lambda\in \Lambda :x^\dagger_\lambda \neq 0 \}$ has sufficiently small cardinality
(that is, it is sufficiently sparse).
Results for $p<1$ can be found in \cite{Cha07,DavGri09,FouLai09,SaaChaYil08}.
\end{remark}

In this section we prove well-posedness of \req{min_prob_sparse}
and derive convergence rates in a possibly infinite dimensional setting.
This inverse problems point of view has so far only been treated
for the case $p=1$ (see~\cite{GraHalSch11}).
The more general setting has only been considered
for Tikhonov regularization
\begin{equation*}
   \norm{\Fo x  - y}^2  + \alpha \rf_{p}(x) \to \min
 \end{equation*}
(see \cite{ChaComPesWaj07,DauDefDem04,Gra09b,GraHalSch08,Lor08,Zar09}).

\subsection{Well-Posedness}
In the following, $\tau$ denotes the weak topology on $\ell^2(\Lambda)$,
and  $\tau_p := \tau_{\rf_p}$ denotes the topology as in Definition \ref{de:tauR}.
Then a sequence $(x_k)_{k \in \N} \subset \ell^2(\Lambda)$ converges to $x \in \ell^2(\Lambda)$
with respect to $\tau_p$ if and only if $x_k \to x$ and $\rf_p(x_k) \to \rf_p(x)$.

\begin{proposition}[Well-Posedness]\label{pr:lp:well_posed}
  Let  $\Fo\colon \ell^2(\Lambda) \to \Y$ be a bounded linear operator with dense range.
  Then constrained $\ell^p$ regularization with $0 < p < 2$ is well-posed:
  \begin{enumerate}
  \item {\em Existence:} %\label{wp-sparse:1}
    For every $\respar > 0$ and $y \in \Y$,
    the set of regularized solutions $\sset(\Fo,y,\respar)$ is non-empty.
  \item {\em Stability:}\label{wp-sparse:2}
    Let $(\respar_k)$ and $(y_k)$ be sequences with $\respar_k \to \respar > 0$ and
    $y_k \to  y \in \Y$. Then $\emptyset \neq \Limsuptr_{k\to\infty} \sset(\Fo,y_k,\respar_k) \subset
    \sset(\Fo,y,\respar)$.
  \item {\em Convergence:}\label{wp-sparse:3}
    Let $\norm{y_k - y} \le \respar_k \to 0$ and assume that
    the equation $\Fo x = y$ has a solution in $\ell^p(\Lambda)$.
    Then we have  \[\emptyset \neq \Limsuptr_{k\to\infty} \sset(\Fo,y_k,\respar_k) \subset \sset(\Fo,y,0) \,. \]
    Moreover, if the equation $\Fo x = y$ has a unique $\rf_p$-minimizing solution $x^\dagger$,
    then we have $\Limsuptp_{k\to\infty} \sset(\Fo,y_k,\respar_k) = \{x^\dagger\}$.
  \end{enumerate}
\end{proposition}

\begin{proof}%[of Proposition \ref{pr:well_posed:convex}]
    In order to prove the existence of minimizers, we apply Theorem~\ref{th:existence} by showing that
    $\Phi_\rf( \Fo,y,\respar, t)$ is compact with respect to the weak topology on $\ell^2(\Lambda)$ for every $t>0$
    and is nonempty for some $t$.
    Because $\Fo$ has dense range, the set
    \[
        \Phi_\rf( \Fo,y,\respar, t) = \set{x\in \ell^2(\Lambda):\rf_p(x) \le t, \norm{\Fo(x) - y}^2 \le \respar}
    \]
    is non-empty for $t$ large enough.

    It remains to show that the sets $\Phi_\rf(\Fo,y,\respar, t)$ are weakly compact on
    $\ell^2(\lambda)$ for every positive $t$.
    The functional $\rf_{p}(x) = \sum_{\lambda\in \Lambda} \abs{x_\lambda}^p$ is weakly lower
    semi-continuous (on $\ell^2(\lambda)$)
    as the sum  of non-negative and weakly continuous functionals (see \cite{EkeTem76}).
    Moreover, the mapping $\Fo$ is weakly continuous,
    and therefore $x \mapsto \norm{\Fo x - y}^2$ is weakly lower semi-continuous, too.
    The estimate $\rf_p(x) \ge \norm{x}_{\ell^2(\Lambda)}^p$ (see~\cite[Equation~(5)]{GraHalSch08})
    shows that $\rf_p$ is weakly coercive. Therefore the sets $\Phi_\rf(\Fo,y,\respar, t)$ are
    weakly compact for all $t>0$, see Proposition~\ref{prop:compact}.

    Taking into account Example~\ref{ex:linear1}, it follows that $\rf_p$, $\sif$, and $\Fo$ satisfy
    Assumption~\ref{as:linear}. Consequently, Items~\ref{wp-sparse:2} and~\ref{wp-sparse:3} follow
    from Proposition~\ref{pr:linear_stab}.
\end{proof}

\begin{remark}
In the case $p>1$, the functional $\rf_p$ is strictly convex, and
therefore the $\rf_p$-minimizing solution $x^\dag$ of $\Fo x = y$ is unique.
Consequently the equality \[\Limsuptp_{k\to\infty} \sset(\Fo,y_k,\respar_k) = \bigl\{x^\dagger\bigr\}\] holds  for every $y$
in the range of the operator $\Fo$.
\end{remark}

\begin{remark}
For the convex case $p\geq 1$, it is shown in \cite[Lemma~2]{GraHalSch08} that the $\tau_p$ convergence of a sequence
$x_k$ already implies $\rf_p(x_k-x) \to 0$.
In particular, the topology $\tau_p$ is stronger than the topology induced by
$\enorm_{\ell^2(\Lambda)}$.
A similar result for $0 < p < 1$ has been derived in~\cite{Gra10}.
\end{remark}

\subsection{Convergence Rates}

In the following, we derive two types of convergence rates results with respect the
$\ell^2$-norm: The convergence rate $\mathcal O\bigl(\respar^{1/2}\bigr)$ (for $p \in (1,2)$), and the convergence  rate
$\mathcal O\bigl(\respar^{\min\{1,1/p\}}\bigr)$ (for every $p \in (0,2)$)
for sparse sequences---here and in the following, $x^\dagger \in \ell^2(\Lambda)$ is called sparse,
if
\[
\supp(x^\dagger) := \set{\lambda \in \Lambda: x^\dagger_\lambda \neq 0}
\]
is finite.
The convergence rates results for constrained $\ell^p$  regularization, derived in this section,
are summarized in Table~\ref{tb:rates}.
\begin{table}[htb]
  \centering
  \begin{tabular}{l |@{\hspace{0.4cm}} l |@{\hspace{0.4cm}} l |@{\hspace{0.4cm}} l  }
    \toprule
    {Rate}       &    {Norm}                & {Premises (besides $\range(\Fo^*) \cap \partial \rf_p \neq \emptyset)$} & Result\\
    \midrule
    $\respar^{1/2}$& $\enorm_{\ell^2} $ & $p\in (1,2)$ & Prop.~\ref{pr:rate_lp} \\
    $\respar^{1/2}$& $\enorm_{\ell^p} $ & $p\in (1,2)$ & Rem.~\ref{re:rate_lp}\\
    $\respar^{1/p}$& $\enorm_{\ell^2} $ & $p\in [1,2)$, sparsity, injectivity on $V$ &
    Prop.~\ref{pr:rate_lp_sparse}\\ \midrule
    $\respar$      & $\enorm_{\ell^2} $ &
    \!\!\!\begin{tabular}{l}
      $p\in (0,1)$, uniqueness of $x^\dagger$,\\
      sparsity, injectivity on $V$
    \end{tabular} & Prop.~\ref{pr:rate_lp_sparse}\\
\bottomrule
\end{tabular}
\caption{Convergence rates for constrained $\ell^p$ regularization.\label{tb:rates}}
\end{table}

For $p \ge 1$, the same type of results (Propositions~\ref{pr:rate_lp}, \ref{pr:rate_lp_sparse})
has also been obtained  for $\ell^p$-Tikhonov regularization in~\cite{GraHalSch08,SchGraGroHalLen09}.
The results for the non-convex case, $p\in (0,1)$, are based on~\cite{Gra10},
where the same rate for non-convex Tikhonov regularization with a--priori
parameter choice has been derived (see also~\cite{Gra10b}).
Similar, but weaker, results have been already been derived in~\cite{BreLor09,Gra09b,Zar09}
in the context of Tikhonov regularization.
In~\cite{Zar09}, the conditions for the convergence rates result for non-convex regularization are
basically the same as in Proposition~\ref{pr:rate_lp_sparse},
but only a rate of order $\mathcal O\bigl(\respar^{1/2}\bigr)$ has been obtained.
In~\cite{BreLor09,Gra09b}, a linear convergence rate $\mathcal{O}(\respar)$ is proven,
but with a considerably stronger range condition:
each standard basis vector $e_\lambda$, $\lambda \in \Lambda$,
has to satisfy $e_\lambda \in \range \Fo^*$.

\begin{proposition}\label{pr:rate_lp}
  Let $1< p < 2$, $x^\dagger = (x^\dagger_\lambda)_{\lambda \in \Lambda}  \in \ell^{2}(\Lambda)$,
  and let $\Fo\colon \ell^2(\Lambda) \to \Y$ be a bounded linear operator.
  Moreover, assume that there exists $\omega \in \Y$ with $\partial \rf_p( x^\dagger ) = \Fo^* \omega$.
  Then the set  $\sset(\Fo,y,\respar) =: \set{x_\respar }$ consists of a single element
  and there exists a constant $d_p>0$ only depending on $p$, such that
  \begin{equation}\label{eq:rate_est_lp}
    \mnorm{x_\respar - x^\dagger}_{\ell^{2}(\Lambda)}^2
    \le
    \frac{d_p \norm{\omega}}{3 + 2\rf_p(x^\dag)}\,
    \bigl( \respar + \mnorm{\Fo x^\dagger - y } \bigr)
  \end{equation}
  for all $\respar>0$ and $y \in \Y$ with $\snorm{ \Fo(x^\dagger)-y } \le \respar$.
\end{proposition}

\begin{proof}
  The assumption $\partial \rf_p(x^\dagger) = \Fo^*\omega$ then implies that
  \req{bregman_est} is satisfied with $W = \X^*$,
  $\gamma_1=0$ and $\gamma_2 = \norm{\omega}$.
  Theorem~\ref{thm:rates} therefore implies the inequality
  \begin{equation}\label{eq:rate_lp:2}
    \sup \set{D_{\partial \rf_p (x^\dagger)}(x_\respar,x^\dagger):x_\respar \in \sset(\Fo,y,\respar) }
    \leq
    \norm{\omega} \bigl(\respar+\mnorm{\Fo x^\dagger - y}\bigr)\;.
  \end{equation}
  From~\cite[Lemma~10]{GraHalSch08} we obtain the inequality
  \begin{equation}\label{eq:rate_lp:3}
    \mnorm{x-x^\dagger}_{\ell^{2}(\Lambda)}^2
    \le
    \frac{d_p}{3 + 2\rf_p(x^\dagger) + \rf_p(x)}
    \, D_{\partial \rf_p (x^\dagger)} \bigl(x,x^\dagger\bigr)
  \end{equation}
  for all $x \in \domain(\rf_p)$.
  Now, \req{rate_est_lp} follows from~\req{rate_lp:2} and~\req{rate_lp:3}.
\end{proof}

\begin{remark} \label{re:rate_lp}
Since $\ell^p(\Lambda)$ is 2-convex (see \cite{LinTza79}) and continuously embedded in $\ell^2(\Lambda)$,
Proposition \ref{pr:rates_norm} provides an alternative estimate for $x_\respar - x^\dagger$ in terms of the
stronger distance $\enorm_{\ell^{p}(\Lambda)}$.
The prefactor in \req{rate_est_banach}, however, is constant, whereas the prefactor in \req{rate_est_lp}
tends to $0$ as $\rf_p(x^\dagger)$ increases. Thus the two estimates are somehow
independent from each other.
\end{remark}

\begin{proposition}[Sparse Regularization]\label{pr:rate_lp_sparse}
  Let $p \in (0,2)$, let $x^\dagger = (x^\dagger_\lambda)_{\lambda \in \Lambda} \in \ell^{2}(\Lambda)$ be sparse,
  and let $\Fo\colon \ell^2(\Lambda) \to \Y$ be bounded linear.
  Assume that one of the following conditions holds:
  \begin{itemize}
  \item We have $p \in (1,2)$, there exists $\omega \in \Y$ with $\partial\rf_p(x^\dagger) = \Fo^*\omega$,
    and $\Fo$ is injective on
    \[
    V = \set{x\in \ell^2(\Lambda):\supp(x) \subset \supp(x^\dagger)}\;.
    \]
  \item We have $p = 1$, there exist $\xi = (\xi_\lambda)_{\lambda \in \Lambda}\in \partial\rf_{1}(x^\dagger)$
    and $\omega \in \Y$ with $\xi = \Fo^* \omega$,
    and $\Fo$ is injective on
    \[
    V =  \set{x\in \ell^2(\Lambda): \supp(x)\subset  \set{\lambda \in \Lambda : \abs{\xi_\lambda} = 1 } } \,.
    \]
  \item We have $p \in (0,1)$, $x^\dagger$ is the unique $\rf_p$-minimizing
  solution of $\Fo x = \Fo x^\dagger$, and $\Fo$ is injective on
  \[
  V = \set{x \in \ell^2(\Lambda) : \supp(x) \subset \supp(x^\dagger)}\;.
  \]
  \end{itemize}
  Then
  \begin{multline*}%\label{eq:rate_sparse}
    \sup\set{ \mnorm{x_\respar - x^\dagger}_{\ell^{2}(\Lambda)} :x_\respar \in \sset(\Fo,y,\respar),\,
    \mnorm{\Fo x^\dagger-y}\le\respar}
    \\ = \mathcal O \left(\respar^{\min\set{1,1/p}} \right) \ \text{ as } \respar \to 0 \;.
  \end{multline*}
\end{proposition}

\begin{proof}
  Assume first that $p \in (1,2)$.
  Define $W := \{w(x) := -c\norm{x-\tilde{x}}^p:\tilde{x} \in \X,\ c > 0\}$.
  Then the functional $\rf_p$ is convex at $x^\dagger$ with respect to $W$.
  Moreover it has been shown in~\cite[Proof of Thm.~14]{GraHalSch08}
  that there exists $w(x) = -c\norm{x-x^\dagger}^p \in \partial_W(x^\dagger) \subset W$ such that
  for some $\eta_1$, $\eta_2 > 0$ the inequality
  \begin{equation}\label{eq:rate_sparse_h1}
  -w(x) = c\mnorm{x-x^\dagger}^p \le \eta_1 \bigl(\rf_p(x)-\rf_p(x^\dagger)\bigr) + \eta_2\mnorm{\Fo(x-x^\dagger)}
  \end{equation}
  holds on $\sset(2\respar,y^\dagger,\Fo)$ for $\respar$ small enough.
  Using Remark~\ref{re:bregman_est_var}, Theorem~\ref{thm:rates} therefore implies the rate
  \[
  \sup\set{D_w(x_\respar,x^\dagger):x_\respar \in \sset(\Fo,y,\respar),\, \mnorm{\Fo x^\dagger-y}\le\respar}
  = \mathcal{O}(\respar)
  \ \text{ as } \respar \to 0\;.
  \]
  The assertion then follows from the fact that the norm on $\ell^2(\Lambda)$
  can be bounded by the Bregman distance $D_w$.

  The proofs for $p = 1$ and $p \in (0,1)$ are similar;
  the required estimate~\req{rate_sparse_h1} has been shown for $p=1$ in~\cite[Proof of Thm.~15]{GraHalSch08}
  and for $p \in (0,1)$ in~\cite[Eq.~(7)]{Gra10}.
\end{proof}

\section{Conclusion}
Due to modeling, computing, and measurement errors, the solution of an ill-posed equation
$\Fo( x ) = y$, even if it exists, typically yields unacceptable results.
The residual method replaces the exact solution by the set
$\sset(\Fo,y,\respar) = \argmin\set{\mathcal R(x) :\mathcal S( \Fo(x), y ) \leq \respar }$, where $\mathcal R$ is
a stabilizing functional and $\mathcal S$ denotes a distance measure between $\Fo(x)$ and $y$.
This paper shows that in a very general  setting $\sset(\Fo,y,\respar)$ is stable
with respect to perturbations of the data $y$
and the operator $\Fo$ (Lemma~\ref{le:stability} and Theorem~\ref{thm:stability:2}),
and the regularized solutions converge to $\rf$-minimizing solutions of $\Fo(x)=y$
as $\respar \to 0$ (Theorem~\ref{thm:convergence}).
In particular the stability issue has hardly  been considered so far in the literature.

In the case where $\Fo$ acts between linear spaces  $\X$ and $\Y$, stability and convergence have been shown
under a list of reasonable properties (see Assumption~\ref{as:linear}). These assumptions are satisfied for
bounded linear operators, but also for a certain class of nonlinear operators (Example~\ref{ex:stab_nonlinear}).
If $\Y$ is reflexive, $\X$ satisfies the Radon--Riesz property, $\Fo$ is a closed linear
operator, and $\mathcal R$ and $\mathcal S$ are given by powers of the norms on
$\X$ and $\Y$, the set $\sset(\Fo,y,\respar)$ consists of a single element $x_\respar$.
This element is shown to converge strongly to the minimal norm solution $x^\dagger$ as $\respar \to 0$.
In this special situation, norm convergence has also been shown in~\cite[Theorem~3.4.1]{IvaVasTan02}.

In Section \ref{se:rates} we have derived quantitative estimates (convergence rates) for the difference between
$x^\dagger$ and minimizers $x_\respar\in \sset(\Fo,y,\respar)$ in terms of a (generalized) Bregman distance.
All these estimates hold
provided   $\sif(\Fo(x^\dagger), y) \le \respar$ and a source inequality introduced in \cite{HofKalPoeSch07} is satisfied.
For linear operators, the required source inequality follows from a source wise representation of a
subgradient of $\rf$ at $x^\dagger$. This carries on the result of \cite{BurOsh04} for constrained regularization.
In the special case that $\X$ is an $r$-convex Banach space with $r \geq 2$ and $\rf$ is the $r$-th power of the norm on $\X$,
we have obtained  convergence rates $\mathcal O(\respar^{1/r})$ with respect to the norm.
The spaces $\X = L^p(\Omega)$ for $p \in (1,2]$ are examples of 2-convex Banach spaces, leading to the rate
$\mathcal O\bigl(\sqrt \respar \bigr)$  in those spaces.

As an application for our rather general results we have investigated sparse $\ell^p$ regularization with
$p \in (0,2)$. We have shown well-posedness in both the convex ($p\geq 1$) and the non-convex case
($p < 1$). In addition, we have studied the reconstruction of sparse sequence.
There we have derived the improved convergence rates $\mathcal O\bigl(\respar^{1/p}\bigr)$
for the convex and $\mathcal O(\respar)$ for the non-convex case.

\section*{Acknowledgement}

This work has been supported by the Austrian Science Fund (FWF),
projects 9203-N12 and project S10505-N20.

\appendix
\section{Auxiliary results}
\label{sec:appendix}

\begin{lemma_appendix}\label{le:Phiincl2}
  Assume that $(y_k)_{k\in\N}$ converges $\sif$-uniformly to $y \in \Y$
  and the mappings $\Fo_k\colon \X \to \Y$ converge   locally $\sif$-uniformly to $\Fo\colon \X \to \Y$.

Then, for every $\respar > 0$, $t > 0$ and $\eps > 0$, there exists
  some $k_0 \in \N$ such that
  \begin{equation}\label{eq:Phiincl2}
  \Phi_\rf(\Fo,y,\respar-\eps,t') \subset \Phi_\rf(\Fo_k,y_k,\respar,t')
  \subset \Phi_\rf(\Fo,y,\respar+\eps,t')
  \end{equation}
  for every   $t' \le t$ and $k \ge k_0$.
\end{lemma_appendix}

\begin{proof}
  Since $y_k \to y$ $\sif$-uniformly and $\Fo_k \to \Fo$ locally $\sif$-uniformly,
  there exists $k_0 \in \N$ such that
  \begin{equation}\label{eq:Phiincl2:1}
    \begin{aligned}
      \abs{\sif(\Fo_k(x),y_k) - \sif(\Fo_k(x),y)} &\le \eps/2\,,\\
      \abs{\sif(\Fo_k(x),y) - \sif(\Fo(x),y)} &\le \eps/2\,,
    \end{aligned}
  \end{equation}
  for all $x \in \X$ with $\rf(x) \le t$ and $k \ge k_0$.

  Now let $t' \leq t$ and let $x \in \Phi_\rf(\Fo,y,\respar-\eps,t')$.
  Then~\req{Phiincl2:1} implies that
  \begin{multline*}
  \abs{\sif(\Fo_k(x),y_k) - \sif(\Fo(x),y)}
  \\ \le \abs{\sif(\Fo_k(x),y_k) - \sif(\Fo_k(x),y)}
      + \abs{\sif(\Fo_k(x),y) - \sif(\Fo(x),y)}
      \le \eps\,,
  \end{multline*}
  and thus
  \[
  \sif(\Fo_k(x),y_k) \le \sif(\Fo(x),y) + \eps \le \respar\,,
  \]
  that is, $x \in \Phi_\rf(\Fo_k, y_k, \respar, t')$,
  which proves the first inclusion in~\req{Phiincl2}.
  The second inclusion is shown in a similar manner.
\end{proof}

The following lemma states that the value of the minimization problem~\req{min_prob}
behaves well as the  parameter $\respar$ decreases.

\begin{lemma_appendix}\label{le:Hlimit}
  Assume that $\Phi_\rf(\Fo,y,\gamma,t)$ is $\tau$-compact for every $\gamma$ and every $t$.
  Then the value $v$ of the constraint optimization problem~\req{min_prob} is right continuous in the first variable,   that is,
  \begin{equation}\label{eq:Hlimit}
    v(\Fo,y,\respar) = \lim_{\eps\to 0+} v(\Fo,y,\respar+\eps) 
    = \sup_{\eps > 0} v(\Fo,y,\respar+\eps)\;.
  \end{equation}
\end{lemma_appendix}

\begin{proof}
  Since $\Phi_\rf(\Fo,y,\respar,t) \subset \Phi_\rf(\Fo,y,\respar+\eps,t)$,
  it follows that $v(\Fo,y,\respar) \ge v(\Fo,y,\respar+\eps)$   for every $\eps > 0$,
  and therefore $v(\Fo,y,\respar) \ge \sup_{\eps > 0} v(\Fo,y,\respar+\eps)$.

  \medskip
  In order to show the converse inequality, let $\delta > 0$.
  Then the definition of $v(\Fo,y,\respar)$ implies
  that $\Phi_\rf\bigl(\Fo,y,\respar,v(\Fo,y,\respar)-\delta\bigr) = \emptyset$.
  Since (cf.~Lemma~\ref{le:Phi_prop})
  \begin{equation}\label{eq:Hlimit:1}
  \emptyset = \Phi_\rf\bigl(\Fo,y,\respar,v(\Fo,y,\respar)-\delta\bigr)
  = \bigcap_{\eps > 0} \Phi_\rf\bigl(\Fo, y,\respar+\eps,v(\Fo,y,\respar)-\delta\bigr)
  \end{equation}
  and the right hand side of~\req{Hlimit:1} is a decreasing family of
  compact sets.  It follows that already
  $\Phi_\rf\bigl(\Fo,y,\respar+\eps_0,v(\Fo,y,\respar)-\delta\bigr) = \emptyset$
  for some $\eps_0 > 0$, and thus
  \[
  \sup_{\eps > 0} v(\Fo,y,\respar+\eps) \geq v(\Fo,y,\respar+\eps_0) \ge v(\Fo,y,\respar)-\delta\;.
  \]
  Since $\delta$ was arbitrary, this shows the assertion.
\end{proof}

\begin{lemma_appendix}\label{le:subseq}
  Let $(\sset_k)_{k \in \N}$ be a sequence of subsets of $\X$. Then
  $U = \Limt_{k\to\infty} \sset_k$, if and only if every subsequence $(\sset_{k_j})_{j\in\N}$
  satisfies 
  \[
  U = \Limsupt_{j\to\infty} \sset_{k_j}\;.
  \]
\end{lemma_appendix}

\begin{proof}
  See~\cite[\S 29.V]{Kur66}.
\end{proof}

\bibliographystyle{plain}

\def\cprime{$'$} \providecommand{\noopsort}[1]{}

\end{document}